\documentclass[12pt]{article}
\usepackage{amsmath}
\usepackage{amssymb}
\usepackage{dsfont}
\usepackage{cite}
\newsymbol \blackbox 1004
\newcommand{\eh}{\hfill}\newlength{\sperr}

\newenvironment{proof}{{\settowidth{\sperr}{\bf\rm
Proof}%
\par\addvspace{0.3cm}\noindent\parbox[t]{1.3\sperr}
{\bf\rm P\eh r\eh o\eh o\eh f\eh }%
}}{\nopagebreak\mbox{}
$\blackbox$\par\addvspace{0.3cm}}

\def\nn{\nonumber}
\def\a{\alpha}
\def\al{\aleph}
\def\b{\beta}
\def\g{\gamma}

\def\de{\delta}

\def\vk{\varkappa}

\def\s{\sigma}
\def\la{\lambda}
\def\om{\omega}

\def\t{\theta}

\def\vp{\varphi}
\def\vt{\vartheta}
\def\ve{\varepsilon}
\def\ze{\zeta}
\def\wh{\widehat}
\def\wt{\widetilde}
\def\ov{\overline}
\def\uv{\underline}
\def\br{\breve}

\def\BC{{\mathbb C}}
\def\BD{{\mathbb D}}
\def\BR{{\mathbb R}}
\def\BN{{\mathbb N}}
\def\clp{{\mathcal P}}
\def\cla{{\mathcal A}}
\def\clb{{\mathcal B}}

\def\cld{{\mathcal D}}
\def\cle{{\mathcal E}}
\def\clf{{\mathcal F}}

\def\clh{{\mathcal H}}

\def\clj{{\mathcal J}}

\def\clm{{\mathcal M}}
\def\cln{{\mathcal N}}
\def\clu{{\mathcal U}}

\def\cld{{\mathcal D}}

\def\mfa{{\mathfrak A}}
\def\ker{{\rm Ker}}

\def\diag{\mathrm{diag}}
\newcommand{\E}{\mathrm{e}}
\newcommand{\I}{\mathrm{i}}
\def\mf{\mathfrak}

\newtheorem{Pa}{Paper}[section]
\newtheorem{Tm}[Pa]{{\bf Theorem}}
\newtheorem{La}[Pa]{{\bf Lemma}}
\newtheorem{Cy}[Pa]{{\bf Corollary}}
\newtheorem{Rk}[Pa]{{\bf Remark}}

\newtheorem{Nn}[Pa]{{\bf Notation}}
\newtheorem{Pn}[Pa]{{\bf Proposition}}

\title{Interpolation Khrushchev-type formulas \\ for structured operators,\\ inequalities and asymptotic relations}

\author{Alexander Sakhnovich\footnote{This research    was supported by the
Austrian Science Fund (FWF) grant  DOI: 10.55776/Y963.}}

\date{}
\begin{document}
\maketitle

\begin{flushleft}
Faculty of Mathematics, University of Vienna, \\
Oskar-Morgenstern-Platz 1, A-1090 Vienna, Austria\\
E-mail: oleksandr.sakhnovych@univie.ac.at
\end{flushleft}

\vspace{0.3em}

\begin{abstract} We show that interpolation results in the $S$-nodes theory
may be considered as Khrushchev-type formulas. If separation of the well-known
Verblunsky (Schur) coefficients occurs in Khrushchev formulas, the separation
of the so the called new Verblunsky-type coefficients occurs in the interpolation formulas 
of the $S$-nodes theory. General asymptotic inequalities (and equalities) for the $S$-nodes,
with application to the block Hankel matrices, are derived using this approach.
Another asymptotic inequality, needed in the proofs and important in itself, is derived
in Appendix B.
\end{abstract}

\vspace{0.3em}


\vspace{0.3em}


\section{Introduction}\label{Intro}
\setcounter{equation}{0}
We show that interpolation results in the $S$-nodes theory
may be considered as Khrushchev-type formulas. If separation of the well-known
Verblunsky (Schur) coefficients occurs in Khrushchev formulas, the separation
of the so the called new Verblunsky-type coefficients occurs in the interpolation formulas 
of the $S$-nodes theory. On the new Verblunsky-type coefficients, 
see our works \cite{ALSJApprTh} as well as \cite{ALS-ConstrAppr}
(and some references therein).
 In this paper, new Verblunsky-type coefficients
and interpolation results for the cases of block Toeplitz and Hankel matrices are discussed
in Sections \ref{Prel} and \ref{Sect3}. A general interpolation result for the $S$-nodes
from \cite{SaL2} is recalled in Section~\ref{Int}. 

We note that interpolation theorems and related results  are essential in the first place in the study of the interconnections between 
matrix measures $d\mu$  (or , equivalently, nondecreasing  matrix functions $\mu(t)$) and corresponding
structured operators $S$ whereas Khrushchev formulas highlight connections between $d\mu$ 
and corresponding orthogonal polynomials (see also Remark \ref{Rktype}).
The interrelations between interpolation theorems 
and  matrix measures $d\mu$ are used in order to derive
the ``entropy" (or Arov-Krein) inequalities   in
Section~\ref{Ine}.

Finally, general asymptotic inequalities (and equalities) for the $S$-nodes
are presented in Theorem~\ref{TmMAs} and Corollary \ref{CyUniq}.
An application to the block Hankel matrices is given in Corollary \ref{CyH}.
Sections \ref{Ine} and \ref{As} are related to the seminal note \cite{ArKr1} by
D.Z. Arov and M.G. Krein. For instance, \cite[Thm.~2]{ArKr1} could be used
in the proof of  Theorem~\ref{TmMAs}. Unfortunately, the promised
proofs of the results of \cite{ArKr1} and, in particular, of a somewhat more complicated Theorem~2  there did not appear
(although some analogs of Theorems~1 and~3 from \cite{ArKr1} have been proved in \cite{ArKr2}).
Therefore,  we use in our proof (of Theorem~\ref{TmMAs}) Proposition~\ref{LIn} derived in Appendix \ref{App}.
This proposition is also of independent interest. The results from Appendices \ref{MB} and \ref{Det} are used
in the proof of Theorem \ref{TmMAs} as well.

Now, let us discuss the  literature on the subject in somewhat greater detail (although  only a very small part of the existing literature could  be mentioned here).
Khrushchev formula (see, e.g., \cite{Khru0} and \cite[Thm. 8.67]{Khru}) is a well-known tool in the theory of orthogonal polynomials on the unit
circle and related problems of analysis. An interesting  matrix version of  Khrushchev formula is given in \cite{CeGr} for the case of block CMV matrices in particular
(see \cite[Thm. 4.3]{CeGr}). There, a matrix function (matrix valued function) depending on Verblunsky (or Schur) coefficients $\{\a_k\}_{k=0}^\infty$ is expressed via a matrix function
determined by the Verblunsky coefficients $\{\a_k\}_{k=0}^{n-1}$ and another matrix function determined by the coefficients $\{\a_k\}_{k=n}^{\infty}$. Thus, one may speak about
Khrushchev formulas as a kind of separation of  Verblunsky coefficients.

At the same time an interesting approach to interpolation and spectral  problems is based on the transfer matrix function in Lev Sakhnovich form \cite{ALS87, ALSJFA, SaSaR, SaL1, SaL2, SaL3} 
(and $S$-nodes theory therein) as well as on the
important works by V.P. Potapov \cite{Am, Pot1}. Using factorisations of transfer matrix functions (in particular, transfer matrix functions corresponding to Toeplitz matrices), important applications were obtained in our papers \cite{ALS87, ALSJFA, ALS-ConstrAppr}. Alternative proofs (of some of these  results) via orthogonal polynomials and generalised Khrushchev formula 
were later presented in the interesting papers \cite{DerSi, FKL} (see also \cite{FKM} and some references therein).

 It was asked  in the famous book on the orthogonal polynomials by B.~Simon \cite{SiOPUC}
how the corresponding works of L. Sakhnovich are related to the theory of orthogonal polynomials on the unit circle (OPUC). In the work \cite{ALSJApprTh} (see also a quite recent
paper  \cite{ALS-ConstrAppr} and references therein), we introduced new Verblunsky-type coefficients closely related to the factorisation of the  transfer matrix function
in \cite{ALS87, ALSJFA, SaSaR, SaL1, SaL2, SaL3}. The work \cite{ALSJApprTh} may be considered as an answer to B. Simon's question.

A discussion on exploring the connections between  transfer matrix function and matrix-valued Khrushchev formulas is contained in   \cite[p. 921]{CeGr}
and one of the aims of the present paper is to highlight these connections. As already mentioned above, we show that certain interpolation results in  L.~Sakhnovich form may be considered as the separation of  our new Verblunsky-type coefficients $($i.e., as Khrushchev-type formulas$)$.
In particular, 
similar to the ``separation" of Verblunsky (or Schur) coefficients $\{\a_k\}_{k=0}^\infty$ in \cite[Thm.~4.3]{CeGr},
we split  Verblunsky-type coefficients $\{\rho_k\}_{k=0}^\infty$ in the Interpolation  theorem \ref{TmNKhr} for Toeplitz matrices.

Turning to the measures and orthogonal polynomials on the real line (OPRL case), we note that a matrix version of Khrushchev formula for the OPRL case was presented  in the interesting paper \cite[Section 7]{GruV}. The corresponding Schur parameters (or coefficients) are discussed in \cite[Section 5]{GruV}. Verblunsky-type coefficients
for Hamburger moment problem were introduced by us in \cite[Section~3]{ALSJApprTh}, and  our interpolation Theorem \ref{Hamburg} may be considered as a Khrushchev-type result for the OPRL case.

{\bf Notations.} As usually, $\BN$ stands for the set of positive integers, $\BR$ stands for the real axis and $\BC$ stands for the complex plane. By $\Re(Z)$ and  $\Im(Z)$
we denote the real and imaginary, respectively, parts of the scalars or square matrices. Here, $(Z+Z^*)/2$ is the real part of the matrix
$Z$, $\I(Z^*-Z)/2$ is the imaginary part of $Z$,  $\I$ stands for the imaginary unit ($\I^2=-1$), and  $Z^*$ is the matrix (or operator) 
adjoint to $Z$. The notation $\BC_+$ ($\BC_-$) denotes the open upper (lower) semiplane $\Im(z)>0$ ($\Im(z)<0$).
The symbol $\BC^{n\times p}$ denotes the set of $n\times p$ matrices with complex-valued entries, and $\BC^n=\BC^{n\times 1}$
is the n-dimensional Hilbert space with the complex inner product. The $p\times p$ identity matrix is denoted by $I_p$ and
a standard identity operator is denoted by $I$.
The symbol $\clb(\clh_1, \clh_2)$ denotes the  set of bounded operators acting
between the Hilbert spaces $\clh_1$ and $\clh_2$, and $\clb(\clh)$ abbreviates $\clb(\clh, \clh)$.
We write that $S \geq 0$ if the scalar product $(Sf,f)$ is nonnegative for the operator $S=S^*$
and we write $S>0$ if this scalar product is always positive for $f\not=0$.
By  $\ov{\lim}$ and $\uv{\lim}$ we denote upper and lower limits, $\de_{ij}$ is Kronecker delta,
and $\ker(A)$ denotes the kernel (null space) of  the operator $A$.
We often write that certain relation holds for some Herglotz matrix function $\vp(z)$ when this relation
holds for $\mu(t)$ from the Herglotz representation \eqref{c17} of $\vp(z)$.
\section{Preliminaries}\label{Prel}
\setcounter{equation}{0}
In this section, we present some results on positive-definite Toeplitz matrices,
discrete Dirac systems and Verblunsky-type coefficients, which may be found in
the papers \cite{ALSJApprTh, ALS-ConstrAppr}. Several of these results appeared already in our
earlier papers  \cite{FKRS08, ALS73, ALSJFA}. We note that several notations from
 \cite{ALSJApprTh, ALS-ConstrAppr} are somewhat changed in this paper (see, e.g.,  Remark \ref{Rk2.1}) for greater convenience.
 \subsection{Toeplitz matrices and operator identities}\label{subD}
 Consider self-adjoint block Toeplitz matrix
\begin{align}\label{c0}
S(n)=S(n)^*=\{s_{j-i}\}_{i,j=1}^n,
\end{align}
where $s_k$ are $p\times p$ blocks. For any such $S(n)$ the following matrix identity
is valid (see \cite{ALS73, ALSJFA} and references therein):
\begin{align}& \label{c1}
AS(n)-S(n)A^*= \I \Pi J \Pi^*;  \quad \Pi=\begin{bmatrix}\Phi_1 & \Phi_2\end{bmatrix},  
\end{align}
where 
\begin{align} & \label{c2}
A= \left\{ a_{j-i}^{\,} \right\}_{i,j=1}^n,
          \quad  a_k  =  \left\{ \begin{array}{lll}
                                  0 \, & \mbox{ for }&
k> 0  \\
\displaystyle{\frac{\I}{ {\, 2 \,}}} \,
                                I_p
                                  & \mbox{ for }& k =
0  \\
                                \, \I \, I_p
                                  & \mbox{ for }& k <
0
                          \end{array} \right.
, \qquad J=\begin{bmatrix} 0 & I_p \\ I_p & 0 \end{bmatrix};
\\ & \label{c3}
\Phi_1 = \left[
\begin{array}{c}
I_{p}  \\ I_p \\ \cdots \\ I_{p}
\end{array}
\right], \quad \Phi_2 =  \left[
\begin{array}{l}
s_0/2  \\ s_0/2 + s_{-1} \\ \cdots \\ s_0/2+ s_{-1} + \ldots +
s_{1-n}
\end{array}
\right] +\I \Phi_1 \nu, \quad \nu = \nu^*;
\\ & \label{c3'}
A=A(n), \quad \Pi=\Pi(n), \quad \Phi_1=\Phi_1(n), \quad \Phi_2=\Phi_2(n),
\end{align}
and $\nu$ is a normalising matrix which need not be fixed at the beginning.

Here, we assume that $S(n)>0$ and so $S(n)$ is invertible as well.
The transfer matrix function $w_A$ in Lev Sakhnovich form \cite{SaL1}  
is given, for the case of the Toeplitz matrix $S(n)$ and the identity \eqref{c1}, by the formula:
\begin{align}& \label{c4}
w_A(n, \la)= I_{2p}-\I J\Pi(n)^*S(n)^{-1}\big(A(n)-\la I_{np}\big)^{-1}\Pi(n).
\end{align} 
Since $S(n)>0$, we have $S(k)>0$ $(1 \leq k \leq n)$, and so all the matrices $S(k)$
are invertible and 
\begin{align}& \label{c8}
t_k:=\begin{bmatrix} 0  & \ldots & 0 & I_p \end{bmatrix}S(k)^{-1}\begin{bmatrix} 0  & \ldots & 0 & I_p \end{bmatrix}^*>0.
\end{align}
Introduce also $p \times p$ matrices $X_k$ and $Y_k$ by the equalities
\begin{align}
 & \label{c7}
\begin{bmatrix} X_k & Y_k  \end{bmatrix}=\begin{bmatrix} 0  & \ldots & 0 & I_p \end{bmatrix}S(k)^{-1}\begin{bmatrix} \Phi_1(k) & \Phi_2(k)  \end{bmatrix}.
\end{align}
According to the general factorisation theorem 
for transfer matrix functions $w_A$ \cite{SaL1} (see also \cite[Thm. 1.16]{SaSaR} and further references therein),
we have a factorisation
\begin{align}& \label{c5}
w_A(n,\la)=w_n(\la)w_{n-1}(\la) \ldots w_1(\la) \quad (n \in \BN), 
\end{align}
where
\begin{align}& \label{c6}
w_k(\la):=I_{2p}-\I \left(\frac{\I}{2}-\la\right)^{-1}J\begin{bmatrix} X_k^* \\ Y_k^* \end{bmatrix} t_k^{-1}\begin{bmatrix} X_k & Y_k  \end{bmatrix} .
\end{align}
 \subsection{Discrete Dirac systems}
The factorisation \eqref{c5}, \eqref{c6} is essential in establishing {\it one to one correspondence $($see \cite[Thm. 2.6]{ALSJApprTh}$)$ between the
discussed above pairs of matrices $\{S(n), \nu\}$  $\,(S(n)>0, \,\, \nu=\nu^*)$,  and discrete
Dirac systems}
\begin{align}& \label{c10}
y_{k+1}(z)=\left(I_{2p} + \I z j C_k\right)y_k(z), \quad j:=\begin{bmatrix}I_p & 0 \\ 0 & -I_p \end{bmatrix},
\end{align}
{\it where $0\leq k<n$ and} 
\begin{align} & \label{c11}
C_k>0, \quad C_k j C_k=j \quad (0 \leq k<n).
\end{align}
This correspondence is given (in one direction) by  the formula
\begin{align}& \label{c12}
C_k:=2K^*\b(k)^*\b(k)K-j, \quad \b(k):=t_{k+1}^{-1/2}\begin{bmatrix} X_{k+1} & Y_{k+1}  \end{bmatrix},
\end{align}
where
\begin{align}& \label{c13}
K:= \frac{1}{\sqrt{2}}\begin{bmatrix}I_p & -I_p \\ I_p & I_p \end{bmatrix}, \quad K^*=K^{-1}, \quad K^*J K=j.
\end{align}

The pair $\{S(n),\nu\}$ is recovered from the Dirac system \eqref{c10}, \eqref{c11} using Taylor series \eqref{c19} 
(see also Remark \ref{Rk1to1}). In order to do this, we need some related interpolation results.
The fundamental solution $W_k(\la)$ of the system \eqref{c10},  \eqref{c11} is normalised by the condition
\begin{align} & \label{c9'} 
W_0(z) \equiv I_{2p}.
\end{align}
Then, $W_k(z)$ is connected with the transfer matrix $w_A(k,z)$ via the formula
\begin{align}
& \label{c9}
W_{k}(z)=(1 -\I z)^{k}K^*w_A\big(k,1/(2z)\big)K \quad (0<k \leq n).
\end{align}
\begin{Rk} \label{Rk2.1} Systems $y_{k+1}(\la)=\big(I_{2p} - (\I /\la) j C_k\big)y_k(\la)$
are considered in \cite{ALSJApprTh, ALS-ConstrAppr} instead of systems \eqref{c10} here, that is, $W_k(z)$ here coincides with $W_k(-1/z)$ in the notations
of \cite{ALSJApprTh, ALS-ConstrAppr}.  
Taking into account this
 change, \eqref{c9} follows from \cite[(2.11)]{ALS-ConstrAppr}. 
 \end{Rk}
 Similar to  \cite[(3.5)]{ALS-ConstrAppr}, we set
\begin{equation} \label{c14}
{\mathfrak A}_n(z)={\mathfrak A}(z)= \{{\mathfrak A}_{ij}(z)
\}_{i,j=1}^2=\left(1-\frac{\I}{2} z\right)^{-n}JjK
W_{n}( -\ov{z}/2)^*K^*jJ,
\end{equation}
and introduce  Weyl functions of the discrete Dirac systems \eqref{c10}, \eqref{c11} by the linear fractional transformations  \cite[(3.7)]{ALS-ConstrAppr}:
\begin{equation}\label{c15}
\varphi (z) = \I \bigl( {\mathfrak A}_{11}(z)R(z)+
{\mathfrak A}_{12}(z)Q(z)\bigr) \bigl( {\mathfrak A}_{21}(z ) R(z)+{\mathfrak A}_{22}(z)Q(z)\bigr)^{-1}.
\end{equation}
Here, $z\in \BC_+$,  ${\mathfrak A}_{ij}(z)$ are $p\times p$ blocks of ${\mathfrak A}(z)$, and the pairs $\{R(z), \, Q(z)\}$  are {\it nonsingular, with property-$J$}, that is,
$R(z)$ and $Q(z)$ are meromorphic $p\times p$ matrix functions in $\BC_+$ satisfying relations
 \begin{align} &  \label{c16}
R(z)^*R(z)+Q(z)^*Q(z)>0, \quad \begin{bmatrix}R(z)^* & Q(z)^* \end{bmatrix}J \begin{bmatrix}R(z) \\ Q(z) \end{bmatrix}\geq 0
 \end{align} 
(excluding, possibly, some isolated points $z\in \BC_+$). 

We note that $\vp$ above is denoted by $w$  and $\om$ in \cite{ALS87} and \cite{ALS-ConstrAppr}, respectively.
Hence, according to \cite{ALS87, ALS-ConstrAppr}, each matrix function $\vp(z)$ of the form \eqref{c15}  belongs to Herglotz class.
Therefore, $\vp(z)$ admits Herglotz representation
\begin{align}& \label{c17}
\vp(z )=\g z +\theta+\int_{-\infty}^{\infty}\frac{1+t z }{(t-z )(1+t^2)}d\mu(t); \\
& \label{c18}
 \g\geq 0, \quad \theta=\theta^*, \quad \int_{-\infty}^{\infty}{(1+t^2)^{-1}}{d\mu(t)}<\infty,
\end{align} 
where $\mu(t)$ is a nondecreasing $p\times p$ matrix function. Furthermore, taking into account \cite[(3.9)]{ALS-ConstrAppr} (where
our $\vp$ from \eqref{c15} here is denoted by $\om$)   and
\cite[Thm. 6.2]{FKRS08}, we obtain the following Taylor series at $\zeta=0$:
 \begin{align} &  \label{c19}
-\I \vp\left( 2\I \frac{1-\zeta}{1+\zeta}\right)=\frac{s_0}{2}+\I \nu+\sum_{k=1}^{\infty}s_{-k}\zeta^k.
 \end{align} 
 \begin{Rk} \label{Rk1to1} Note that the $p\times p$ matrices $s_{-k}$ $(0\leq k<n)$ in Taylor series \eqref{c19} are the corresponding
 blocks of the block Toeplitz matrix $S(n)$. On the other hand, according to \cite[Thm. 6.2]{FKRS08}, the coefficients $s_{-k}$ $(k\geq~n)$
 are such that all the matrices $S(N)=\{s_{j-i}\}_{i,j=1}^N$, where $s_k=s_{-k}^*$,  $N>~n$, are nonnegative
 $($i.e., $S(N)\geq 0)$. Moreover, each sequence $\{s_{-k}\}_{k\geq n}$, which determines an extension of $S(n)$ with this property,  is generated by some
$($nonsingular, with property-$J)$ pair $\{R(z), \, Q(z)\}$  and by Taylor series \eqref{c19} of the corresponding matrix function
$-\I \vp$ given by \eqref{c15}.
 \end{Rk}

 \begin{Rk} \label{RkN1to1}
The one to one correspondence between Dirac systems \eqref{c10}, \eqref{c11} and Verblunsky-type $p\times p$ matrix coefficients 
 $\rho_k$:
 \begin{align} &  \label{c20}
\|\rho_k\|<1 \quad (0\leq k<n)
 \end{align} 
  is simpler $($see  \cite{ALSJApprTh}$)$. It is given by the so following representations  of $C_k$ $($so called Halmos extensions of $\rho_k) :$
  \begin{align}\label{c21}&
C_k= { \mathcal D_k}  F_k, \quad {\mathcal D}= {\mathrm{diag}}\Big\{
\big(
I_{p}-  \rho_k  \rho_k^* \big)^{-\frac{1}{2}}, \, \,
\big(I_{p}- \rho_k^* \rho_k\big)^{-\frac{1}{2}}\Big\}, 
\\ \label{c22}&
 F_k= \left[
\begin{array}{cc}
I_{p} &  \rho_k \\   \rho_k^* & I_{p}
\end{array}
\right] \quad (0\leq k<n),
\end{align}
and by the formula
\begin{align}\label{c23}&
\rho_k=\left(\begin{bmatrix} I_p & 0 \end{bmatrix}C_k\begin{bmatrix} I_p \\ 0 \end{bmatrix}\right)^{-1}\begin{bmatrix} I_p & 0 \end{bmatrix}C_k
\begin{bmatrix} 0 \\ I_p \end{bmatrix} \qquad (0\leq k<n),
\end{align}
which easily follows from \eqref{c21} and \eqref{c22}.
\end{Rk}

\section{Khrushchev-type formulas}\label{Sect3}
\setcounter{equation}{0}
\subsection{Khrushchev-type formulas \\ and block Toeplitz matrices}
{\bf 1.} Let us assume that an infinite sequence  of $p\times p$ matrices $s_{-k}$ $(0 \leq k< \infty)$ is given. We also
assume that all the matrices $S(N)=\{s_{j-i}\}_{i,j=1}^N$, where $s_k=s_{-k}^*$,  $N \geq 1$, are positive definite
 $($i.e., $S(N)> 0)$. Recall that slightly more general sequences (where $S(N)\geq 0$) were discussed in Remark \ref{Rk1to1}.
 It follows from \cite[p.~194]{ALSJApprTh} that, in the case $N>n>0$, Dirac system corresponding to $S(N)$ is an extension  (to the interval $0\leq k<N$)
 of the Dirac system corresponding to $S(n)$. In other words, the first $n$ coefficients $C_k$ (and $\rho_k$) for both systems coincide.
 \begin{Nn}\label{NnCln}
 The family of Weyl functions given by \eqref{c15} is denoted by $\cln\big({\mf A}_n\big)$.
 \end{Nn}
 According to \cite[pp. 651, 652]{ALS-ConstrAppr}, the families $\cln\big({\mf A}_n\big)$ are decreasing, that is, 
 $\cln\big({\mf A}_n\big) \supset \cln\big({\mf A}_N\big)$ for $N>n>0$. Moreover, there is a unique Weyl function $\vp_{\infty}(z)$
 for Dirac system \eqref{c10} on the semiaxis $0 \leq k<\infty$ (or, equivalently, for our infinite sequence   $s_{-k}$ $(0 \leq k< \infty)$):
 \begin{align}\label{c24}&
\bigcap_{n\geq 1}\cln\big({\mf A}_n\big)=\{\vp_{\infty}(z)\}.
\end{align} 
\begin{Rk}\label{vpinfty}
One can see that $\vp_{\infty}(z)$ is determined by our sequence  of $p\times p$ matrices $s_{-k}$ $(0 \leq k< \infty)$ $($and by the matrix $\nu=\nu^*$) via
formula \eqref{c19} as well as by by the  sequence  of $p\times p$ Verblunsky-type coefficients $\rho_{k}$ $(0 \leq k< \infty)$
via formula \eqref{c24}, where the matrix functions $\mfa_n$ are given by \eqref{c14} and \eqref{c10}, \eqref{c21}, \eqref{c22}.
 In this way,  a one to one correspondence  between an infinite sequence
$s_{-k}$ $(0 \leq k< \infty)$ such that all $S(N)=\{s_{j-i}\}_{i,j=1}^N>0$ $($and $\nu=\nu^*)$ on one side and  an infinite sequence  $\rho_{k}$ $(0 \leq k< \infty)$, where $\|\rho_k\|<1$,
on the other side, is established $($see \cite[Thm.~2.12]{ALSJApprTh}$)$.
\end{Rk}
{\bf 2.} Let  us split the given sequence of Verblunsky-type coefficients $\rho_{k}$ $\break (0 \leq k< \infty)$ into two sequences:
\begin{align}\label{c25}&
\{\rho_k\}_{k=0}^\infty=\{\rho_k\}_{k=0}^{n-1}\cup \{\wt \rho_k\}_{k=0}^\infty \quad (\wt \rho_k:=\rho_{k+n}).
 \end{align} 
\begin{Nn}\label{Nntvp}
By $\wt \vp_{\infty}(z)$, we denote Weyl function of Dirac system \eqref{c10} determined by the Verblunsky-type coefficients $\wt \rho_k$ $(0\leq k<\infty)$.
\end{Nn}
Now, an analog of  S. Khrushchev's result  (of Khrushchev formula) may be formulated.
\begin{Tm}\label{TmNKhr} Let a sequence $\{\rho_k\}_{k=0}^\infty$ of the $p\times p$ matrix coefficients $\rho_k$ such that $\|\rho_k\|<1$ $(0\leq k<\infty)$ be given.
Let ${\mf A}_n= \{{\mathfrak A}_{ij}(z)
\}_{i,j=1}^2$ be given by the formula \eqref{c14}, where $W_n(z)$ is the fundamental solution $($at $k=n)$ of the Dirac system \eqref{c10}, \eqref{c11}
determined by the Verblunsky-type coefficients $\{\rho_k\}_{k=0}^{n-1}$. Then, we have
\begin{equation}\label{c26}
\varphi_{\infty} (z) = \I \bigl( {\mathfrak A}_{11}(z)\big(-\I \wt \varphi_{\infty} (z)\big)+
{\mathfrak A}_{12}(z)\bigr) \bigl( {\mathfrak A}_{21}(z ) \big(-\I \wt \varphi_{\infty} (z)\big)+{\mathfrak A}_{22}(z)\bigr)^{-1},
\end{equation}
where $\varphi_{\infty}$ is determined by the Verblunsky-type coefficients $\{\rho_k\}_{k=0}^\infty$ and $\wt \varphi_{\infty}$ is determined by the Verblunsky-type coefficients 
$\{\wt \rho_k\}_{k=0}^\infty$ $(\wt \rho_k:=\rho_{k+n})$.
\end{Tm}
\begin{proof}. By $\wt W_N(z)$, we denote  the fundamental solution $($at $k=N$, normalised by $\wt W_0(z)\equiv I_{2p})$ of the Dirac system determined by the 
Verblunsky-type coefficients $\wt \rho_k$ $(0\leq k<\infty)$. The matrix function $\wt {\mf A}_N(z)$ is obtained by the substitution of $\left(1-\frac{\I}{2} z\right)^{-N}
\wt W_{N}( -\ov{z}/2)^*$ instead of 
$\left(1-\frac{\I}{2} z\right)^{-n}
W_{n}( -\ov{z}/2)^*$ into the right-hand side of \eqref{c14} (and by the substitution of $\wt {\mf A}_N(z)$ instead of 
$ {\mf A}_n(z)$ into the left-hand side of \eqref{c14}). We will add natural numbers standing in the indices of the so called frames (i.e., matrices of coefficients of our linear
fractional transformations)
$ {\mf A}_n(z)$ and $\wt {\mf A}_N(z)$ as variables in the notations of their $p\times p$ blocks and will write, for instance,
${\mathfrak A}_n(z)= \{{\mathfrak A}_{ij}(n,z)
\}_{i,j=1}^2$ and $\wt {\mathfrak A}_N(z)= \{\wt {\mathfrak A}_{ij}(N, z)
\}_{i,j=1}^2$.

It follows from \eqref{c24} (for the case of Verblunsky-type coefficients $\wt \rho_k$) that
$\wt \vp_{\infty}(z)\in \cln\big(\wt {\mf A}_N\big)$ for any $N\in \BN$. In other words (see Notation \ref{NnCln}), we have
\begin{align}\nn
-\I \wt \varphi_{\infty} (z) = &\bigl( \wt{\mathfrak A}_{11}(N,z)\wt R(N,z)+
\wt {\mathfrak A}_{12}(N,z)\wt Q(N,z)\bigr) 
\\ & \label{c27} \times
\bigl( \wt {\mathfrak A}_{21}(N, z ) \wt R(N,z)+\wt {\mathfrak A}_{22}(N,z)\wt Q(N, z)\bigr)^{-1}
\end{align}
for some nonsingular, with property-$J$ pairs $\{\wt R(N,z), \, \wt Q(N, z)\}$. Formula \eqref{c27}
may be rewritten in the form
\begin{align}\label{c28}&
\begin{bmatrix}
-\I  \wt \varphi_{\infty} (z)\\ I_p
\end{bmatrix}=\wt {\mf A}_N(z) \begin{bmatrix}
\wt R(N,z) \\ \wt Q(N,z)
\end{bmatrix}\bigl( \wt {\mathfrak A}_{21}(N, z ) \wt R(N,z)+\wt {\mathfrak A}_{22}(N,z)\wt Q(N, z)\bigr)^{-1}.
 \end{align} 
Let us denote the right-hand side of \eqref{c26}
by $\psi(z)$. Then, \eqref{c28} yields
\begin{align}\nn
\begin{bmatrix}
-\I \psi(z) \\ I_p
\end{bmatrix}= & {\mf A}_n(z) \wt {\mf A}_N(z) \begin{bmatrix}
\wt R(N,z) \\ \wt Q(N,z)
\end{bmatrix}
\\ & \nn \times
\bigl( \wt {\mathfrak A}_{21}(N, z ) \wt R(N,z)+\wt {\mathfrak A}_{22}(N,z)\wt Q(N, z)\bigr)^{-1}
\\ \label{c29}& \times 
\bigl( {\mathfrak A}_{21}(z ) \big(-\I \wt \varphi_{\infty} (z)\big)+{\mathfrak A}_{22}(z)\bigr)^{-1}.
 \end{align} 
 Since $\wt \rho_k=\rho_{k+n}$, one can see that $\wt W_N(z)W_n(z)=W_{n+N}(z)$.
 Hence, formula \eqref{c14} implies
 \begin{align}& \label{c30}
 {\mf A}_n(z) \wt {\mf A}_N(z)={\mf A}_{n+N}(z).
\end{align}
It follows from \eqref{c29}, \eqref{c30} (and Notation \ref{NnCln}) that $\psi(z) \in \cln\big({\mf A}_{n+N}\big)$.
That is, $\psi(z) \in \cln\big({\mf A}_{n+N}\big)$ for any $N\in \BN$. Thus, in view of \eqref{c24} we have
$\break \psi(z)=\varphi_{\infty} (z)$, where $\psi(z)$ denotes the right-hand side of \eqref{c26}.
\end{proof}
\begin{Rk}\label{RkGen}  Theorem \ref{TmNKhr} and related to it considerations in Introduction mean that 
the solutions of  interpolation problems of  the form \eqref{c15}, where $\mf A$ is expressed via the
transfer matrix function $w_A$ $($see, e.g., \cite{ALS87, ALSJApprTh, SaL2} and references therein$)$
may be considered as generalisations of the Khrushchev-type formulas corresponding to the
new Verblunsky-type coefficients.
\end{Rk}
\subsection{Khrushchev-type formulas \\ and block Hankel matrices}\label{Hank}
Formula \eqref{c19} holds in the unit disk $|\zeta|< 1$,  and  the matrix version of Khrushchev formula is considered in the OPUC case \cite{CeGr} 
in this disk as well. 

In the OPRL case, the authors of  \cite{GruV} consider Herglotz (Nevanlinna)
matrix functions in $\BC_+$. They fix some measure $d\mu$ on the real line and corresponding Jacobi matrices
 \begin{align}& \label{c31}
\clj_n=\begin{bmatrix}
b_0 & a_0 & 0 & \ldots & 0 \\
a_0^* & b_1 & a_1  & 0 &\ldots \\
0 & \ddots & \ddots & \ddots & \ddots \\
0 & \ldots & a_{n-2}^* & b_{n-1} &  a_{n-1} \\
0 & \ldots & 0 & a_{n-1}^* & b_n
\end{bmatrix},
\end{align}
where $a_k$ and $b_k$ are $p\times p$ matrices, $a_k=a_k^*>0$, $b_k=b_k^*$.
The coefficients 
\begin{align}& \label{c31'}
b_0,\, a_0, \, b_1, \, a_1, \, b_2, \, a_2, \, \ldots
\end{align}
are called (in \cite{GruV}) Schur (Verblunsky) parameters of the corresponding
Herglotz matrix function $\vp(z)$ with some  matrix  measure $d\mu$ in Herglotz representation \eqref{c17} of $\vp$. The block Hankel matrices
\begin{align}& \label{c32}
H(n)=\{H_{i+j-2}\}_{i,j=1}^n \quad (n \geq 1),
\end{align}
where 
\begin{align}& \label{c33}
H_k=\int_{-\infty}^{\infty}t^{k}d\mu(t) ,
\end{align}
don't appear in \cite{GruV} although the requirement  
\begin{align}& \label{c34}
\int_{-\infty}^{\infty}t^{2n-2}d\mu(t) <\infty 
\end{align}
is given.
\begin{Rk}\label{Rktype}
We note that Verblunsky coefficients and Khrushchev formulas are
closely connected with the orthogonal polynomials corresponding to the measure $d\mu$
and are somewhat less connected with the structured matrices $($Toeplitz, Hankel, etc.$)$ generated by this measure.
At the same time, transfer matrix functions, Verblunsky-type coefficients and Khrushchev-type
formulas are directly connected with the structured matrices.
\end{Rk}
Let us recall  a solution of the interpolation problem  (that is, generalised Krushchev-type result in the spirit of Remark  \ref{RkGen})
for the case of a block Hankel matrix $H(n)$ (i.e., matrix of the form \eqref{c32} where $H_k$ are $p\times p$ blocks).
 We refer to \cite{SaL2} (or to a more convenient  for our purposes presentation in \cite{ALSJApprTh}).
 
 The self-adjoint Hankel matrix $H=H(n)=H(n)^*$ (equivalently, the Hankel matrix where $H_k=H_k^*$) satisfies the matrix identity
 \begin{align}& \label{H2}
AH-HA^*=\I \Pi J \Pi^*, \quad \Pi=\Pi(n)=\begin{bmatrix}\Phi_1(n) & \Phi_2(n)\end{bmatrix},
\end{align}\
where $J$ is given in \eqref{c2} and
 \begin{align}
& \label{H3}
A=A(n)=\{a_{ij}\}_{i,j=1}^n, \quad a_{ij}=\de_{i-1,j}I_p, \\ & \label{H4}
\Phi_1=\Phi_1(n)=-\I \begin{bmatrix} 0 \\ H_0 \\ H_1 \\ \cdots \\ H_{n-1}\end{bmatrix}, \quad 
\Phi_2=\Phi_2(n)=\begin{bmatrix} I_p \\  0\\  0\\ \cdots \\ 0\end{bmatrix}.
\end{align}
Clearly $A(n)$ and $H(n)$ are $n p \times np$ matrices and $\Phi_1(n)$ and $\Phi_2(n)$
are $ np \times p$ matrices.  In view of \eqref{H2}, an explicit solution of the interpolation problem  (i.e., of the truncated
Hamburger moment problem) for $H>0$ is easily obtained via the method of operator identities \cite{SaL1, SaL2, SaL3} 
(see also the references therein)
using V.P. Potapov's fundamental matrix inequalities
\cite{Am, Pot1}. For that purpose we introduce the transfer matrix function
\begin{align}& \label{H5}
w_A(n, \la)= I_{2p}-\I J\Pi(n)^*H(n)^{-1}\big(A(n)-\la I_{np}\big)^{-1}\Pi(n),
\end{align} 
 and the frame ${\mf A}_n$ with the $p\times p$ blocks $\mfa_{ij}$:
\begin{align}& \label{H7}
\mfa_n(z)=\{\mfa_{ij}(z)\}_{i,j=1}^2=w_A(n,1/\ov{z})^*.
\end{align}
Now, we have the following theorem.
 \begin{Tm}\label{Hamburg} Assume that the block Hankel matrix $H=H(n)$ is positive-definite
$($i.e., $H>0)$. 

Then, the matrix functions $\vp(z)$ given by the linear fractional
transformations
\begin{align}& \label{H8}
\varphi (z) = \I \bigl( {\mathfrak A}_{11}(z)R(z)+
{\mathfrak A}_{12}(z)Q(z)\bigr) \bigl( {\mathfrak A}_{21}(z ) R(z)+{\mathfrak A}_{22}(z)Q(z)\bigr)^{-1},
\end{align}
where $\{R(z),\, Q(z)\}$ are nonsingular, with property-$J$ pairs $($see \eqref{c16}$)$, belong to the Herglotz class, that is, $\Im \big(\vp(z)\big) \geq 0$ for $z\in \BC_+$.
Moreover, each matrix function $\vp(z)$ admits a unique Herglotz representation of the form
\begin{align}& \label{H10}
\vp(z)=\int_{-\infty}^{\infty}(t-z)^{-1}d\mu(t) <\infty,
\end{align}
and the $p\times p$ matrix measure $d\mu$ $($or, equivalently, nondecreasing matrix function $\mu(t))$ in this representation satisfies \eqref{c34}.

These and only these  matrix functions $\mu(t)$  $($i.e., $\mu(t)$ given by \eqref{H8}, \eqref{H10}$)$ satisfy the equalities
\begin{align}& \label{H11}
H_k=\int_{-\infty}^{\infty}t^{k}d\mu(t) \quad (0 \leq k <2n-2),
\end{align}
and the inequality 
\begin{align}& \label{H12}
H_{2n-2} \geq \int_{-\infty}^{\infty}t^{2n-2}d\mu(t). \end{align}
\end{Tm}
Relations \eqref{H10} and \eqref{H11} yield a corollary below (see \cite[p. 199]{ALSJApprTh}).
 \begin{Cy} For  $z$ tending non-tangentially to infinity in $\BC_+$, we have an asymptotic expansion $($compare with  \eqref{c19} for Toeplitz matrices$):$
\begin{align} & \label{H12+}
\vp(z)=-\sum_{k=0}^{2n-3}\frac{1}{z^{k+1}}H_k+ O\left(\frac{1}{z^{2n-1}}\right).
\end{align}
\end{Cy}
Let us recall Verblunsky-type coefficients corresponding to $H(n)>0$ (see \cite[Subsection 3.2]{ALSJApprTh}). Similar to the case of the Toeplitz matrices, the transfer matrix
function \eqref{H5} corresponding to Hankel matrix $H(n)$ admits factorisation of the form \eqref{c5}:
\begin{align}& \label{H13-}
w_A(n,\la)=w_n(\la)w_{n-1}(\la) \ldots w_1(\la) \quad (n \in \BN).
\end{align}
Here (for $0 \leq k<n$), we have
\begin{align}& \label{H13}
 w_{k+1}( \la)=I_{2p}+\frac{\I}{\la}J Q_{k}, \quad Q_k= \om_k^* t_{k+1}^{-1}\om_k, 
\end{align}
where $\om_k$ is a $p\times 2p$ matrix and $t_{k+1}$ is a $p\times p$ matrix:
\begin{align}& \label{H14}
 \om_k: =P_2(k+1) T(k+1)\Pi(k+1), \quad T(r):=H(r)^{-1},
 \\ & \label{H15}
P_2(r):=
\begin{bmatrix} 0 & \ldots & 0 & I_p
 \end{bmatrix} \in \BC^{p\times pr},
 \quad t_r:=P_2(r)T(r)P_2(r)^*>0.
\end{align}
In view of \eqref{H13-}, discrete canonical system
\begin{align}& \label{H16}
y_{k+1}(\la)=w_{k+1}(\la)y_k(\la) \quad (0 \leq k<n)
\end{align}
corresponds to $H(n)$. The matrices $\om_k$ have the properties
\begin{align}
 & \label{H17}
 \om_kJ \om_k^*=0, \quad \I \om_{k}J \om_{k-1}^*=t_{k+1} \quad (0 < k <n);  \quad \om_0=\begin{bmatrix} 0 & t_1 \end{bmatrix}.
\end{align}
\begin{Rk}\label{VerH} The matrices $\om_k$ are called  Verblunsky-type coefficients corresponding to $H(n)>0$ \cite{ALSJApprTh}.
In view of \eqref{H13} and \eqref{H17}, they determine the factors $w_k(\la)$, transfer matrix function $w_A(n,\la)$ and canonical system \eqref{H16}.
Moreover, they uniquely determine $H(n)>0$. 

More precisely, according to \cite[Thm. 3.9]{ALSJApprTh}, 
each sequence $\{\om_k\}_{k=0}^{n-1}$ of $2p\times p$ matrices $\om_k$, such that
$$ \om_kJ \om_k^*=0, \quad \I \om_{k}J \om_{k-1}^*>0 \quad (0 < k <n);  \quad \om_0=\begin{bmatrix} 0 & t \end{bmatrix} \quad (t>0),$$
is a sequence of Verblunsky-type coefficients for some Hankel matrix $\break H(n)>0$, and there is a one to one correspondence between the sequences of 
Verblunsky-type coefficients $\{\om_k\}_{k=0}^{n-1}$ and Hankel matrices $H(n)>0$.
\end{Rk}
We note that  Verblunsky coefficients \eqref{c31'} introduced in  \cite{GruV} may be presented as $p\times 2p$
block matrices $\begin{bmatrix} b_k & a_k \end{bmatrix}$, where $b_k=b_k^*$ and $a_k>0$ instead of the conditions \eqref{H17} for our Verblunsky-type
coefficients.

\section{Interpolation formulas for $S$-nodes} \label{Int}
\setcounter{equation}{0}
A much more general (than in the previous sections) interpolation result is formulated in \cite{SaL2} in terms of the so called $S$-nodes (see \cite{SaL1,SaL2, SaL3} and some references
therein). That is, we consider some finite or infinite-dimensional Hilbert space $\clh$ and operators $A,S\in \clb(\clh)$, $\Pi \in \clb(\BC^{2p},\clh)$, which satisfy
the operator identity
\begin{align}& \label{S1}
AS-SA^*=\I \Pi J \Pi^* \quad (S=S^*),
\end{align}
where $J$ is given in \eqref{c2}. {\it Such a triple $\{A,S,\Pi\}$ forms an $S$-node.} (In fact, it is a so called symmetric $S$-node, but
only  symmetric $S$-nodes are studied in this paper.) We partition $\Pi$ into two blocks
\begin{align}& \label{S2}
\Pi=\begin{bmatrix}\Phi_1 & \Phi_2\end{bmatrix}, \quad \Phi_k\in \clb(\BC^{p},\clh) \quad (k=1,2).
\end{align}
\begin{Rk}\label{RkS} For the examples of operator identities \eqref{S1} see, for instance, \eqref{c1} or \eqref{H2}. If one wants
to consider a class of structured operators, one fixes $A$ and $\Phi_2$ $($or $\Phi_1$ in \eqref{c3}$)$, but $S$ and $\Phi_1$
$($$\Phi_2$ in \eqref{c3}$)$ vary. Further we assume that  $A$ and $\Phi_2$ are fixed.
\end{Rk}
\begin{Nn}\label{NnA}
We assume that $\s(A)$ is a finite or countable set of points and denote this class of operators $A$ by $\clf$ $($i.e., $A\in \clf)$.
\end{Nn}
Let $\mu(t)$ be a nondecreasing $p\times p$ matrix function on $\BR$ and $\g, \theta$ be $p\times p$ matrices, where $\g \geq 0, \quad \theta=\theta^*$.
Set 
\begin{align}& \label{S3}
\wt S=S_{\mu}+FF^*, \quad S_{\mu}:=\int_{-\infty}^{\infty}(I-tA)^{-1}\Phi_2[d\mu(t)]\Phi_2^*(I-tA^*)^{-1},
\end{align}
where $F$ is given by the equality
\begin{align}& \label{S4}
AF=\Phi_2\g^{1/2}.
\end{align}
\begin{Nn}  A nondecreasing  matrix function $\mu(t)$ belongs to the class $\cle$ if the integral in \eqref{S3} weakly converges and
\begin{align}& \label{S5}
\int_{-\infty}^{\infty}\frac{d \mu(t)}{1+t^2}<\infty .
\end{align}
\end{Nn}
If  $\mu \in \cle$, the integral below weakly converges as well, and we introduce $\wt \Phi_1$:
\begin{align}& \label{S6}
\wt \Phi_1:= \Phi_{1,\mu}+\I(\Phi_2 \theta+F\g^{1/2}), 
\\ & \label{S7} \Phi_{1,\mu}:=-\I \int_{-\infty}^{\infty}\left(A(I-tA)^{-1}+ \frac{t}{1+t^2}I\right)\Phi_2 d\mu(t).
\end{align}
{\it Here $($see \cite{SaL2}$)$, the interpolation problem is the problem of describing the triples $\{\g\geq 0,\, \t=\t^*,\, \mu \in \cle\}$ such that}
\begin{align}& \label{S8}
S=\wt S, \quad \Phi_1=\wt \Phi_1 .
\end{align}
In order to present the solution of the interpolation problem we need formulas \eqref{H5}, \eqref{H7} in the general
$S$-node setting. Namely, we introduce the frame ${\mf A}(S,z)$ with the $p\times p$ blocks $\mfa_{ij}(S,z)=\mfa_{ij}(z)$:
\begin{align}& \label{S9}
\mfa(S,z)=\{\mfa_{ij}(z)\}_{i,j=1}^2=w_A(1/\ov{z})^*,  \quad w_A(\la)= I_{2p}-\I J\Pi^*S^{-1}\big(A-\la I\big)^{-1}\Pi,
\end{align}
where $w_A$ is the transfer matrix function of the corresponding $S$-node.  Sometimes, we will write
$\mfa_{ij}(S,z)$ instead of $\mfa_{ij}(z)$ (meaning the blocks of $\mfa(S,z)$). We will require that
\begin{align}& \label{S9+}
S \geq \ve I \quad {\mathrm{for\,\, some}} \quad \ve >0, \quad \ker \, \Phi_2=0.
\end{align}
Then, Theorem 1.4.2 and Proposition 1.3.2
from \cite{SaL2}
yield the following interpolation theorem.
\begin{Tm}\label{TmIntS} Let an $S$-node $\{A,S,\Pi\}$ be given, where $A\in \clf$ and zero is not an eigenvalue of $A$.
Assume that \eqref{S9+} holds. Then, the set of solutions $\{\g, \, \theta, \, \mu(t)\}$
of the interpolation problem \eqref{S8} coincides with the set of $\{\g, \, \theta, \, \mu(t)\}$ in the Herglotz representations \eqref{c17}
of linear fractional transformations
\begin{equation}\label{S10}
\varphi (z) = \I \bigl( {\mathfrak A}_{11}(z)R(z)+
{\mathfrak A}_{12}(z)Q(z)\bigr) \bigl( {\mathfrak A}_{21}(z ) R(z)+{\mathfrak A}_{22}(z)Q(z)\bigr)^{-1}.
\end{equation}
where the pairs $\{R(z), \, Q(z)\}$  are nonsingular, with property-$J$.
\end{Tm}
\begin{Rk}\label{RkIntH} In the case of Hankel matrices, zero is an eigenvalue of $A$. Thus, in this case
one needs some modification of Theorem \ref{TmIntS}. Indeed, we have an inequality $($instead of the equality$)$ in
\eqref{H12}.
\end{Rk}
Formula \eqref{S10} is similar to \eqref{c15} and, correspondingly,  the set of $\vp$ given by
\eqref{S10} is denoted by $\cln\big(\mfa(S)\big)$.
\section{Inequalities   for the $S$-nodes}\label{Ine}
\setcounter{equation}{0}
 It easily follows from \eqref{S9} and identity \eqref{S1} (see also \cite{SaL1} or \cite[(1.88)]{SaSaR}) that
\begin{align}& \label{S12}
\mfa(S,z)J \mfa(S,\ov{\la})^*=J-\I (z-\la)\Pi^*(I-zA^*)^{-1}S^{-1}(I-\la A)^{-1}\Pi.
\end{align} 
Interpolation theorem \ref{TmIntS} and this formula are essential in the study of the interconnections between 
matrix measures $d\mu$  (or , equivalently, nondecreasing  matrix functions $\mu(t)$) and corresponding
structured operators $S$.

In particular, an important characteristic of $S$ and $S^{-1}$ is the matrix function
\begin{align}
&  \label{I8} 
\rho(z,\ov{z}):=
\I (\ov{z}-z)\Phi_2^*(I-zA^*)^{-1}S^{-1}(I-\ov{z} A)^{-1}\Phi_2 .
\end{align} 
Relations \eqref{S12} and \eqref{I8} yield 
\begin{align}
&  \label{I8+} 
\rho(z,\ov{z})= \mfa_{21}(z)\mfa_{22}(z)^*+\mfa_{22}(z)\mfa_{21}(z)^*.
\end{align} 
Under condition \eqref{S9+}, we have $\rho(z,\ov{z})>0$ (in the points $z\in \BC_+$ of the invertibility of $I-zA^*$).
We note the matrix function $\rho(z,\ov{z})$ is not related to the Verblunsky-type coefficients $\rho_k$
(using $\rho$ in both cases, we preserve the notations from the important here works \cite{ALS87, ALSarxFact, ALSJApprTh}).
\begin{Rk}\label{RkLF}
Relations \eqref{I8+} and $\rho(z,\ov{z})>0$ imply $($see \cite[Proposition~1.43]{SaSaR}$)$ that the matrix function
${\mathfrak A}_{21}(z ) R(z)+{\mathfrak A}_{22}(z)Q(z)$ is invertible in $\BC_+$ $($excluding, possibly, isolated points$)$.
Thus, taking into account \eqref{S12}, we see that linear fractional transformations \eqref{S10} are well defined and generate
Herglotz functions if only \eqref{S9+} holds and $\mfa(S,z)$ given by \eqref{S9} is meromorphic in $\BC_+$.
\end{Rk}
We will need some results  from \cite{ALSarxFact}. 
Recall that  functions $\break z=(\overline{z_0}\zeta -z_0 )(\zeta -1)^{-1}$ $\,\,(z_0\in \BC_+)$
bijectively map the unit disk $\BD$ $(|z|<1)$ onto $\BC_+$. The choice of $z_0\in \BC_+$ is not essential for our  purposes and we choose $z_0=\I$.
For some  function or matrix function $G(z)$ $(z\in \BC_+)$ the accent  ``widehat'' means (in the notations of \cite{ALSarxFact}, which we also use here) the mapping into the corresponding
function on $\BD$, that is,
\begin{align}& \label{A1}
\wh G(\ze)=G\big(\I(1+ \zeta  )\big/(1-\zeta)\big) \quad (\ze \in \BD).
\end{align}
Let us recall the necessary definitions, notations and properties
for the functions holomorphic in the unit disk $\BD$ (see \cite{ALSarxFact} or \cite{KiKa, Priv} and 
other references therein).  An outer (or maximal in the terminology of \cite{Priv}) function $h(\ze)$   is an analytic in $\BD$ function,
which admits representation
\begin{align}& \label{A2!} 
h(\ze)=\E^{\I \eta} \exp\left\{\frac{1}{2\pi}\int_{0}^{2\pi}\ln\big( g(\vt)\big)\frac{\E^{\I \vt}+ \ze}{\E^{\I \vt}- \ze} d\vt \right\} \quad (g(\vt)\geq 0),
\end{align}
where $\eta \in \BR$ and $\ln\big( g(\vt)\big)$ is integrable on $[0,2\pi]$. 
Clearly, $1/h$ and the product of outer  functions are outer as well. Representation \eqref{A2!} yields that $ g(\vt)=\big|h\big(\E^{\I \vt}\big)\big|$.
A holomorphic in $\BD$ function $f$ belongs to Smirnov class $D$ if it may be represented as a ratio of a function from the Hardy class $H^{\infty}(\BD)$ and of an 
outer function from $H^{\infty}(\BD)$. According to \cite[Thm. 4.29]{RoRov} (or \cite[Lemma 2.1]{KiKa}) the function $h$ is outer if and only if $h,1/h \in D$.
The notation $D^{(p\times p)}$ stands for the  class of $p\times p$ matrix functions with the entries belonging to  $D$.
A matrix function $f(\ze)$ is called an outer matrix function 
if  $f\in D^{(p\times p)}$ and its determinant  is an outer function. A $p\times p$ matrix function $f(\ze)$ is outer if and only if
$f,\, f^{-1}\in D^{(p\times p)}$.
\begin{Tm}\label{TmArKr} \cite[Thm. 3.1]{ALSarxFact} Let an $S$-node $\{A,S, \Pi\}$, such that 
the operators $I-zA^*$ have  bounded inverses for $z\in \BC_+$ $($excluding, possibly, some isolated points, which are the poles of the
matrix function $\mfa(S,z)$ introduced in \eqref{S9}$)$, be given. Let relations \eqref{S9+} and 
\begin{align}& \label{B19}
\wh \mfa_{21}(\ze), \, \wh \mfa_{21}(\ze)^{-1} \in D^{(p\times p)}
\end{align}
hold. Assume that $\vp(z)\in \cln(\mfa(S))$ and the
matrix function $\mu(t)$ in Herglotz representation \eqref{c17}, \eqref{c18} of $\vp$ satisfies Szeg{\H o} condition 
\begin{align}& \label{S11}
\int_{-\infty}^{\infty}(1+t^2)^{-1}\ln\big(\det\mu^{\prime}(t)\big)dt>-\infty,
\end{align} 
where  $\mu^{\prime}$ is the positive semi-definite derivative of the absolutely continuous part of $\mu$.

Then, $\mu^{\prime}$ admits factorisation $($unique up to a constant unitary factor from the left$)$:
\begin{align}& \label{A4}
\mu^{\prime}(t)=G_{\mu}(t)^*G_{\mu}(t),
\end{align} 
where $G_{\mu}(t)$ is the boundary value function of  $G_{\mu}(z )$ $(z\in \BC_+)$, the entries of $\wh G_{\mu}(\ze )$ belong 
to the Hardy class $H^2(\BD)$ $($thus, also belong to $D)$ and $\,\det\big(\wh G_{\mu}(\ze )\big)$
is an outer function.
For this $G_{\mu}(z )$, we have
\begin{align}& \label{B13!}
2 \pi  G_{\mu}(z)^* G_{\mu}(z)\leq \rho\big(z, \ov{z}\big)^{-1}  \quad (z \in \BC_+).
\end{align} 
The equality in \eqref{B13!} holds at some point $z=\la \in \BC_+$ if and only if our $\vp(z)\in \cln(\mfa(S))$ is generated $($see \eqref{S10}$)$ by the
constant pair $\{R,Q\}:$
\begin{align}& \label{B31}
R(z)\equiv \mfa_{22}(S,\la)^*, \quad Q(z)\equiv \mfa_{21}(S,\la)^*.
\end{align} 
\end{Tm}
Sufficient conditions for \eqref{B19} to hold are presented in the next proposition.
\begin{Pn}\label{PnPM}   \cite[Proposition 4.1]{ALSarxFact} Let an $S$-node $\{A,S,\Pi\}$ be given  and let
the operators $I-zA$ have bounded inverses for $z$ in the domains  $\{z: \, \Im(z)\leq 0\}$ and $\{z: \, \Im(z)>0, \,\, |z| \geq r_0\}$ for some $r_0>0$.
Assume that  relations \eqref{S9+}  hold. Finally, let
\begin{align}& \label{Z2}
{\ov{\lim}}_{r\to \infty}\big(\ln(\clm(r))\big/r^{\vk}\big)<\infty  
\end{align} 
for some  $0<\vk<1$ and $\clm(r)$ given by
\begin{align}& \label{Z3}
\clm(r)=\sup_{r_0<|z|<r} \|(I-z A)^{-1}\|.
\end{align} 
Then, we have \eqref{B19}, that is, $\wh \mfa_{21}(\ze)$ is an outer matrix function.
\end{Pn} 
From Theorems \ref{TmIntS} and \ref{TmArKr}, and from Proposition \ref{PnPM} follows the corollary
below.
\begin{Cy}\label{CyIneq} Let an $S$-node $\{A,S,\Pi\}$ be given and let
the operators $I-zA$ have bounded inverses for $z$ in the domain   $\{z: \, \Im(z)>0, \,\, |z| \geq r_0\}$ for some $r_0>0$
as well as for $z$ in the domain  $\{z: \, \Im(z)\leq 0\}$.
Assume that zero is not an eigenvalue of $A$ and that relations \eqref{S9+}   and \eqref{Z2} hold.
Let $\{\g, \vt, \mu\}$ be a solution of the interpolation problem \eqref{S8} and let $\mu^{\prime}(t)$
satisfy  Szeg{\H o} condition \eqref{S11}. Then, the inequality \eqref{B13!} is valid.
\end{Cy}
\section{Asymptotic properties   for the $S$-nodes}\label{As}
\setcounter{equation}{0}
Let $\clh_1 \subset \clh_2 \cdots \subset \clh_r \cdots$ be a given sequence of imbedded Hilbert spaces $\clh_k$ $(1\leq k<\infty)$ with complex inner products.
By $P_k$, we denote the orthogonal projectors from  $\clh_r$ $(r>k)$ onto $\clh_k$ and, by $P_k^{\perp}$, we denote 
the orthogonal projectors from  $\clh_r$ onto the   orthogonal complement of $\clh_k$ in $\clh_r$.
In this section, we assume that a sequence of the $S$-nodes $\{A_k,S_k,\Pi_k\}$
 is given such that  $A_k,S_k\in \clb(\clh_k);\,$ $\, \Pi_k=\begin{bmatrix}\Phi_1(k) & \Phi_2(k)\end{bmatrix}$, $\, \Phi_i(k) \in\clb(\BC^{p}, \clh_k)$ and 
\begin{equation} \label{As1}
A_k=P_kA_rP_k^*, \,\, S_k=P_kS_rP_k^*, \,\, \Pi_k=P_k\Pi_r, \,\,P_kA_r\big(P_k^{\perp}\big)^*=0 \,\, {\mathrm{for}}\,\, r>k.
\end{equation} 
\begin{Pn}\label{PnNond} Let the operators $S_k$ satisfy the inequalities $S_k\geq \ve_k I$ $\break (\ve_k>0)$.
Then, the sequences of matrices
\begin{align}
&  \label{I8!} 
\rho_k(z,\ov{z}):=
\I (\ov{z}-z)\Phi_2(k)^*(I-zA_k^*)^{-1}S_k^{-1}(I-\ov{z} A_k)^{-1}\Phi_2(k), 
\end{align} 
where $z\in \BC_+$ and $I-zA_k^*$ have bounded inverses,
are nondecreasing.
\end{Pn}
Clearly, formula \eqref{I8} for $\rho$ will coincide with \eqref{I8!} if we apply \eqref{I8} to the $S$-node $\{A_k,S_k,\Pi_k\}$.
\begin{proof} of Proposition \ref{PnNond}.  Factorisation theorem for the transfer matrix function \cite{SaL1} have been already  used in Sections \ref{Prel} and \ref{Sect3}
for symmetric $S$-nodes in the cases of Toeplitz and Hankel matrices. In the case of the $S$-node $\{A_r,S_r,\Pi_r\}$ $(r>k)$
and corresponding transfer matrix function $w_{A_r}(\la)$, we have
\begin{align}& \label{As2}
w_{A_r}(\la)=\breve w_A(\la)w_{A_k}(\la),
\end{align} 
where $\breve w_A(\la)$ is the transfer matrix function corresponding to the $S$-node $\{A_{22},T_{22}^{-1}, T_{22}^{-1}P_k^{\perp}\Gamma\}$ and
\begin{align}& \label{As3}
A_{22}:=P_k^{\perp}A_r\big(P_k^{\perp}\big)^*, \quad T_{22}:=P_k^{\perp}(S_r)^{-1}\big(P_k^{\perp}\big)^*, \quad \Gamma:=S_r^{-1}\Pi_r.
\end{align} 
We note that the inequality  $T_{22}, T_{22}^{-1} \geq \breve \ve I$ holds (under proposition's conditions) for some $\breve \ve >0$.
Recall that, according to \eqref{S9}, we have $\break \mfa(S,z)=w_A\big(1/ \ov{z}\big)^*$, and set  $\breve \mfa(S_r,z):=\breve w_A\big(1/ \ov{z}\big)^*$.
Then, formula \eqref{As2} yields
\begin{align}& \label{As4}
\mfa(S_r,z)=\mfa(S_k,z)\breve \mfa(S_r,z).
\end{align} 
Taking into account \eqref{S12} and the inequality  $T_{22}^{-1} \geq \breve \ve I$, we obtain
\begin{align}& \label{As5}
 \breve \mfa(S_r,z)J \breve \mfa(S_r,z)^*\geq J.
\end{align} 
Relations \eqref{As4} and \eqref{As5} imply that
\begin{align}& \label{As6}
 \mfa(S_r,z)J \mfa(S_r,z)^*\geq   \mfa(S_k,z)J \mfa(S_k,z)^*.
 \end{align} 
 In view of \eqref{I8+} and \eqref{As6}, we see that $\rho_r(z,\ov{z})\geq \rho_k(z,\ov{z})$ for $r>k$
 (in the points of invertibility of $I-z A_r^*$ in $\BC_+$).
\end{proof}
Together with $\rho_k(z,\ov{z})$, the matrix functions $\rho_k(\ov{z}, z)$ ($z\in \BC_+$) are also of interest (see, e.g., formula \eqref{B0}).
Similar to \eqref{I8+} we obtain
\begin{align}
\nn 
\rho_k(\ov{z},z)&=
\I (z-\ov{z})\Phi_2(k)^*(I-\ov{z}A_k^*)^{-1}S_k^{-1}(I-{z} A_k)^{-1}\Phi_2(k)
\\   \label{R1} &=
\mfa_{21}(S_k,\ov{z})\mfa_{22}(S_k,\ov{z})^*+\mfa_{22}(S_k,\ov{z})\mfa_{21}(S_k,\ov{z})^*<0
\end{align} 
for $z\in \BC_+$. It follows  from \eqref{S12} that 
\begin{align}& \label{As5!}
 \mfa(S_k,\ov{z})J  \mfa(S_k,\ov{z})^*\leq J,\quad
  \breve \mfa(S_r,\ov{z})J \breve \mfa(S_r,\ov{z})^*\leq J \quad {\mathrm{for}} \,\, z \in \BC_+.
\end{align} 
Therefore, the factorisation formula \eqref{As4} yields the inequality 
\begin{align}& \label{As6!}
 \mfa(S_r,\ov{z})J \mfa(S_r,\ov{z})^*\leq   \mfa(S_k,\ov{z})J \mfa(S_k,\ov{z})^* \quad (r>k,\,\, z\in\BC_+).
 \end{align}
\begin{Cy}\label{CyOvz}
Let the operators $S_k$ satisfy the inequalities $S_k\geq \ve_k I$ $(\ve_k>0)$ and let the operators $I-zA_k$ $(z\in \BC_+)$
have bounded inverses at $z$. 

Then, we have
\begin{align} \label{R2} &
0<-\rho_k(\ov{z},z)\leq -\rho_r(\ov{z},z) \quad (k\leq r, \,\, z \in \BC_+).
\end{align} 
\end{Cy}
\begin{Rk}\label{RkEmb} Relations \eqref{As4} and \eqref{As5} also imply  that 
\begin{align}& \label{As6+}
\cln\big( \mfa(S_r)\big)\subseteq \cln\big(\mfa(S_k)\big) \quad {\mathrm{for}} \quad r>k.
 \end{align} 
\end{Rk}
Now, we formulate and prove our main asymptotical theorem.
\begin{Tm}\label{TmMAs} $I.$ Let the conditions of Proposition \ref{PnPM} be fulfilled for the $S$-nodes
$\{A_k,S_k,\Pi_k\}$ $(1 \leq k < \infty)$. Assume that
\begin{align}& \label{As7}
\vp(z)\in \bigcap\limits_{k\geq 1} \cln\big(\mfa(S_k)\big),
 \end{align} 
and the
matrix function $\mu(t)$ in Herglotz representation \eqref{c17}, \eqref{c18} of $\vp$ satisfies Szeg{\H o} condition \eqref{S11}. Then,
\begin{align}& \label{As8}
\lim\limits_{k\to \infty}\rho_k\big(z, \ov{z}\big)^{-1} \geq 2 \pi  G_{\mu}(z)^* G_{\mu}(z) \quad (z \in \BC_+).
\end{align} 
For all $\la\in \BC_+$, we also have
\begin{align}& \label{As8+}
\lim\limits_{k\to \infty}\det\big(\rho_k\big(\la, \ov{\la}\big)^{-1}\big) >0 .
\end{align} 
$II.$ Let the conditions of Proposition \ref{PnPM} be fulfilled for the $S$-nodes
$\{A_k,S_k,\Pi_k\}$ $(1 \leq k < \infty)$, let
$\mfa(S_k,z)$ be meromorphic in $\BC_-$ and assume that the
inequality \eqref{As8+} holds for some
fixed $\la\in \BC_+$.

Then, there is $\vp(z)$ satisfying \eqref{As7} and Szeg{\H o} condition \eqref{S11}, such that the
 equality holds in \eqref{As8} for this $\vp(z)$ and this $\la\, :$
\begin{align}& \label{As9}
\lim\limits_{k\to \infty}\rho_k\big(\la, \ov{\la}\big)^{-1} = 2 \pi  G_{\mu}(\la)^* G_{\mu}(\la) \quad (\la \in \BC_+).
\end{align} 
\end{Tm}
\begin{proof}. {\it Step 1.} Since $\ker \, \Phi_2(k)=0$, we have $\rho_k\big(z, \ov{z}\big)>0$ for $\rho_k$ given by
\eqref{I8!} and $z\in \BC_+$. Thus, the matrix functions $\rho_k\big(z, \ov{z}\big)$ in \eqref{As8} are invertible (and $\rho_k\big(z, \ov{z}\big)^{-1}>0$).
Now, the existence of the limit on the left-hand side of \eqref{As8} follows from Proposition \ref{PnNond}.
Finally, Theorem \ref{TmArKr} and Proposition~\ref{PnPM} imply the inequality \eqref{As8}. Since $\wh G_{\mu}(\ze)$
is an invertible (outer) matrix function, \eqref{As8+} follows (for all $\la\in \BC_+$) from \eqref{As8}.

 {\it Step 2.}  In order to prove \eqref{As9} (under  theorem's conditions), we consider the matrix functions $\vp_k(z)$
generated (via \eqref{S10}) by the so called ``frames" $\mfa(S_k,z)$ and pairs
\begin{align}& \label{As10}
R_k(z)\equiv R_k= \mfa_{22}(S_k,\la)^*, \quad Q_k(z)\equiv Q_k= \mfa_{21}(S_k,\la)^* 
\end{align} 
for an arbitrary fixed value of $\la \in \BC_+$. According to \eqref{S9}, we have
\begin{align}& \label{As11}
\mfa(S_k,z)=I_{2p}-\I z\Pi_k^*(I-zA_k^*)^{-1}S_k^{-1}\Pi_k J.
\end{align}
Thus, the   frames $\mfa(S_k,z)$   are holomorphic in $\BC_+\cup \BR$ and  the pairs $\{R_k,Q_k\}$ are well defined.
Moreover, taking into account \eqref{I8+} and \eqref{As10}, we obtain a strict inequality
\begin{align}& \label{As12}
R_k^*Q_k+Q_k^*R_k=\rho_k(\la,\ov{\la})>0.
\end{align}
We also have 
\begin{align}& \label{As12+}
\mfa(S_k,z)J \mfa(S_k,z)^*=\mfa(S_k,z)^*J\mfa(S_k,z)=J \quad {\mathrm{for}} \quad z\in \BR.
\end{align}
In view of \eqref{As12} and \eqref{As12+}, the matrix function
\begin{align}& \label{As13-}
F_k(z):={\mathfrak A}_{21}(S_k, z ) R_k+{\mathfrak A}_{22}(S_k,z)Q_k
\end{align}
is invertible on the real axis (see \cite[Proposition 1.43]{SaSaR}). Thus, $\vp_k(z)$ of the form \eqref{S10} may be continuously
extended to the real axis, and, taking again into account \eqref{As12} and  \eqref{As12+}, we obtain
\begin{align}
\nn
\vp_k(z)-\vp_k(z)^*&=\I \big(F_k(z)^{-1}\big)^*\begin{bmatrix}R_k^* & Q_k^*\end{bmatrix}\left(\begin{bmatrix}\mfa_{21}(z)^* \\ \mfa_{22}(z)^*\end{bmatrix}
\begin{bmatrix}\mfa_{11}(z) & \mfa_{12}(z)\end{bmatrix}
\right.
\\ &\nn \quad \left.
+\begin{bmatrix}\mfa_{11}(z)^* \\ \mfa_{12}(z)^*\end{bmatrix}
\begin{bmatrix}\mfa_{21}(z) & \mfa_{22}(z)\end{bmatrix}\right)\begin{bmatrix}R_k \\ Q_k\end{bmatrix}F_k(z)^{-1}
\\ \nn &
=\I \big(F_k(z)^{-1}\big)^*\begin{bmatrix}R_k^* & Q_k^*\end{bmatrix}\mfa(S_k,z)^*J\mfa(S_k,z)\begin{bmatrix}R_k \\ Q_k\end{bmatrix}F_k(z)^{-1}
\\  \label{As13} &
=\I \big(F_k(z)^{-1}\big)^* \rho_k(\la,\ov{\la})F_k(z)^{-1} \quad {\mathrm{for}} \quad z\in \BR.
\end{align}
By $\mu_k$, we denote here the matrix function $\mu$ in the Herglotz representation \eqref{c17} of $\vp_k$. It follows from the well-known Stieltjes-Perron inverse formula 
and from \eqref{As13} that  
$\mu_k(t)$ is absolutely continuous and that 
 \begin{align}& \label{As14}
\mu_k^{\prime}(t)=\frac{1}{2\pi \I}(\vp_k(t)-\vp_k(t)^*)=\frac{1}{2\pi }\big(F_k(t)^{-1}\big)^* \rho_k(\la,\ov{\la})F_k(t)^{-1}.
\end{align}
On the other hand, denoting $\mu_k^{\prime}$ by $\clp_k$ and using \eqref{c17} we easily derive for $\eta>1$ that
 \begin{align}& \label{As15}
\frac{1}{2\I}\eta (\vp_k(\I\eta)-\vp_k(\I\eta)^*)=\eta^2 \g_k +\int_{-\infty}^{\infty}\frac{\eta^2 \clp_k(t)}{t^2+\eta^2}dt\geq 
\int_{-\infty}^{\infty}\frac{ \clp_k(t)dt}{1+t^2}dt.
\end{align}
According to \eqref{As6+}, we have (for all $k$) that $\vp_k(z) \in \cln\big(\mfa(S_1,z)\big)$.
Hence, Proposition \ref{PnMB} yields the uniform boundedness of $\vp_k(\I \eta)$ and so of the left-hand side
of \eqref{As15} (at least for the fixed values $\eta \geq r_0$).
It follows that the right-hand side of \eqref{As15} is uniformly bounded as well and conditions \eqref{Ap1}
of Proposition~\ref{LIn}
are satisfied for $\clp_k(t):=\mu_k^{\prime}(t)$.

 {\it Step 3.}  Taking into account  Remark \ref{HellyM}, we see that there is a convergence \eqref{Ap2} for some subsequence $\clp_{k_r}(t)$ $(r \in \BN)$
and some $\g_{\clp}$ and $\mu$. For simplicity, we will write sometimes $\vp_r(z)$, $\clp_{r}(t)$ and so on instead of  $\vp_{k_r}(z)$, $\clp_{k_r}(t)$...
Proposition \ref{LIn} implies inequality \eqref{Ap3} for the subsequence $\clp_{r}(t)$.

The mentioned above uniform boundedness of $\vp_k(\I\eta)$ yields that the real and imaginary parts of $\vp_k(\I\eta)$,
that is, the expressions
\begin{align}&\label{As16}
(\vp_k(\I\eta)+\vp_k(\I\eta)^*)/2=\t_k+(1-\eta^2)\int_{-\infty}^{\infty}\frac{td\mu_k(t)}{(t^2+\eta^2)(1+t^2)},
\\ &\label{As17}
(\vp_k(\I\eta)-\vp_k(\I\eta)^*)/(2 \I)=\eta \g_k+\int_{-\infty}^{\infty}\frac{\eta d\mu_k(t)}{t^2+\eta^2},
\end{align}
are uniformly bounded as well. In particular, $\g_r$ and $\t_r$ for the mentioned above subsequence $\vp_r(z)$
are uniformly bounded. Hence, we may choose this subsequence in such a way that the limits below exist  and we have
\begin{align}& \label{As18}
\lim_{r\to \infty}\g_r=\g_{\infty} \quad (\g_{\infty}\geq 0), \quad \lim_{r\to \infty}\t_r=\t =\t^*.
\end{align}
It follows from \eqref{As18} and \eqref{Ap2} that there is a limit
\begin{align}& \label{As19}
\lim_{r\to \infty}\vp_r(z)=\vp_{\infty}(z):=\g z +\t+\int_{-\infty}^{\infty}\frac{1+t z }{(t-z )(1+t^2)}d\mu(t), 
\end{align}
where $\g:=\g_{\infty}+\g_{\clp}$. 

Since the matrix functions $\mfa(S_r,z)$ are meromorphic in $\BC_-$, it follows from \eqref{B1} that 
$\mfa(S_r,z)^{-1}$ are meromorphic in $\BC_+$.  Therefore, according to
 \eqref{As19},
$\vp_{\infty}(z)$ is a Herglotz matrix function holomorphic in $\BC_+$, it belongs to the matrix balls (disks)
described in Proposition \ref{PnMB} and is
generated by the frames $\mfa(S_r,z)$ and the corresponding meromorphic pairs 
$$
\begin{bmatrix}\breve R_r(z) \\ \breve Q_r(z)\end{bmatrix}=\mfa(S_r,z)^{-1}\begin{bmatrix}-\I \vp_{\infty}(z) \\ I_p\end{bmatrix}.$$
Clearly, the pairs $\{\br R_r, \br Q_r\}$ are nonsingular, and property-$J$ for them is immediate from \eqref{B4}.
Thus, we have $\vp_{\infty}(z)\in\cln\big(\mfa(S_r)\big)$ (for all $r$ in our subsequence). In view of the  relations \eqref{As6+},
$\vp(z)=\vp_{\infty}(z)$ satisfies \eqref{As7} (for the whole sequence of the embedded sets $\cln\big(\mfa(S_k)\big)$).

 {\it Step 4.}  Recalling now that the inequality \eqref{Ap3} holds for $\clp_{r}(t)=\mu_r^{\prime}(t)$, where $\mu(t)$ in \eqref{Ap3} belongs to the
Herglotz representation of $\vp_{\infty}(z)$, and setting $f(t)=\Im(\la)(1+t^2)|t-\la|^{-2}$ in \eqref{Ap3},  we obtain
\begin{align}&\nn
\int_{-\infty}^{\infty}
 \Im(\la)\ln\Big(\det\big(\mu^{\prime}(t)\big)\Big)\frac{dt}{|t-\la|^2}
 \\ & \label{As21}
 \geq
\displaystyle{\ov{\lim}_{r\to \infty}}\int_{-\infty}^{\infty} \Im(\la)\ln\Big(\det\big(\clp_r(t)\big)\Big)\frac{dt}{|t-\la|^2}. \end{align} 
Let us consider the right-hand side of \eqref{As21}. It follows from \eqref{As14} that
\begin{align}& \label{As28}
\det\big(\clp_r(t)\big)=\det \big(\mu_r^{\prime}(t)\big)=\frac{\det\big(\rho(\la,\ov{\la})\big)}{(2\pi)^p\big|\det\big(F_r(t)\big)\big|^2}.
\end{align}
According to \eqref{As10} and \eqref{As12}, $\mfa_{21}(S_r,z)$ $(z\in \BC_+)$ and $Q_r$ are invertible and, moreover,
$\Re\big(R_rQ_r^{-1}+\mfa_{21}(S_r,z)^{-1}\mfa_{22}(S_r,z)\big)>0$. Hence, we also have $\Re\big(R_rQ_r^{-1}+\mfa_{21}(S_r,z)^{-1}\mfa_{22}(S_r,z)\big)^{-1}>0$.
Therefore, it follows from Smirnov's theorem (see, e.g., \cite[p. 93]{Priv}) that 
$$ h^*\big(R_rQ_r^{-1}+\wh \mfa_{21}(S_r,\ze)^{-1}\wh \mfa_{22}(S_r,\ze)\big)^{\pm 1}h\in H^{\delta}$$ 
for any $h\in \BC^p$ and  the Hardy classes $H^{\delta}$ $(0<\delta <1)$. (Recall that the functions $\wh \mfa_{ij}$ above
are introduced using \eqref{A1}.) Since $H^{\delta} \subset D$, we obtain
\begin{align}& \label{As29}
\big(R_rQ_r^{-1}+\wh \mfa_{21}(S_r,\ze)^{-1}\wh \mfa_{22}(S_r,\ze)\big)^{\pm 1}\in D^{p\times p},
\end{align}
and so $\det \big(R_rQ_r^{-1}+\wh \mfa_{21}(S_r,\ze)^{-1}\wh \mfa_{22}(S_r,\ze)\big)$ is an outer function.
Furthermore, taking into account Proposition \ref{PnPM} and \eqref{As10}, we see that  $\wh \mfa_{21}(S_r,\ze)$ is an outer matrix function
and $Q_r$ is a constant matrix. Thus, for $F_r$ given by \eqref{As13-} we derive
\begin{align}& \label{As30}
\wh F_r(\ze)^{\pm 1}=\big(\wh \mfa_{21}(S_r,\ze)\big(R_rQ_r^{-1}+\wh \mfa_{21}(S_r,\ze)^{-1}\wh \mfa_{22}(S_r,\ze)\big)Q_r\big)^{\pm 1}\in D^{p\times p}.
\end{align}
Using \eqref{A1}, we have 
\begin{align}& \label{As22-}
F_r(\la)=\wh F_r(\ze_0), 
\end{align}
where
$\det \wh F_r(\ze)$ is an outer function and $ \la=\I(1+ \zeta_0  )\big/(1-\zeta_0)$, that is,
\begin{align}& \label{As22}
\zeta_0=(\la -\I)\big/ (\la+\I).
\end{align}
Setting  $h(\ze)=\det \wh F_r(\ze)$ in \eqref{A2!}, we obtain
\begin{equation} \label{As23} 
\ln |\det \wh F_r(\ze_0)|= \frac{1}{2\pi}\int_{0}^{2\pi}\ln\big(\big|\det \wh F_r(\ze)\big| \big)\Re\left(\frac{\ze+ \ze_0}{\ze- \ze_0} \right)d\vt 
\quad (\ze=\E^{\I \vt}).
\end{equation}
Next, we switch in \eqref{As23} from $\ze=\E^{\I \vt}$ to $t=\I(1+ \zeta  )\big/(1-\zeta)$. Similar to \eqref{As22}, we derive
$\ze=(t-\I)\big/ (t+\I)$. It follows that
\begin{align}&\nn
\ze -\ze_0=(t-\I)\big/ (t+\I)-(\la -\I)\big/ (\la+\I)=2\I(t-\la)(\la+\I)^{-1}(t+\I)^{-1},
\\  & \label{As24}
|\ze -\ze_0|^{-2}=(1/4)(\la +\I)(\ov{\la}-\I)(1+t^2)|t-\la|^{-2}.
\end{align}
Thus, we have
\begin{align} \nn
\Re\left(\frac{\ze+ \ze_0}{\ze- \ze_0} \right)&=(1/2)\big((\ze+ \ze_0)(\ov{\ze}- \ov{\ze_0})+(\ov{\ze}+ \ov{\ze_0})(\ze- \ze_0)\big)|\ze -\ze_0|^{-2}
\\ & \nn
=(1/4)\big(1-|\ze_0|^2\big)(\la +\I)(\ov{\la}-\I)(1+t^2)|t-\la|^{-2}
\\ &  \nn
=(1/4)\big((\la +\I)(\ov{\la}-\I)-(\la -\I)(\ov{\la}+\I)\big)(1+t^2)|t-\la|^{-2}
\\ & \label{As25}
= \Im(\la)(1+t^2)|t-\la|^{-2}.
\end{align}
It is easy to see (and follows also from \cite[(2.6)]{ALSarxFact}) that
\begin{align}& \label{As26}
\frac{dt}{1+t^2}=(1/2)d\vt.
\end{align}
Relations \eqref{As22-}, \eqref{As23}, \eqref{As25} and \eqref{As26} imply that
\begin{align}& \label{As31}
\ln |\det F_r(\la)|=\frac{1}{2\pi}\int_{-\infty}^{\infty}
 \Im(\la)\ln\Big( |\det F_r(t)|^2\Big)\frac{dt}{|t-\la|^2}.
\end{align}
Formulas \eqref{As10}, \eqref{As12} and \eqref{As13-} yield
\begin{align}& \label{As32}
F_r(\la)=\rho_r(\la,\ov{\la}).
\end{align}
It is easy to see that 
\begin{align}& \label{As33}
\int_{-\infty}^{\infty}
 \frac{ \Im(\la)dt}{|t-\la|^2}=\pi.
\end{align}
It follows from \eqref{As28} and \eqref{As31}--\eqref{As33} that
\begin{align} \nn
\int_{-\infty}^{\infty} \Im(\la)\ln\Big(\det\big(\clp_r(t)\big)\Big)\frac{dt}{|t-\la|^2}
&= -2\pi \ln\big( \det \big(\rho_r(\la,\ov{\la})\big)\big)-\pi \ln(2\pi)^{p}
\\ & \nn
\quad +\pi \ln\big( \det\big(\rho_r(\la,\ov{\la})\big)\big)
\\ & \label{As34}
=\pi \ln \big(\det\big(2\pi \rho_r(\la,\ov{\la})\big)^{-1}\big).
\end{align} 
{\it Step 5.}  In view of \eqref{As21} and \eqref{As34}, formula \eqref{As8+} implies
that Szeg{\H o} condition \eqref{S11} is fulfilled, indeed, for our $\vp=\vp_{\infty}$ and
the corresponding $\mu$.

Next, let us consider $\ln |\det G_{\mu}(\la)|$ for $\mu$ from the Herglotz representation \eqref{As19}.
According to Theorem \ref{TmArKr}, $\det\big(\wh G_{\mu}(\ze)\big)$ is an outer function. Thus,
 similar to \eqref{As31} we obtain 
 \begin{align}& \label{As27}
\ln |\det G_{\mu}(\la)|=\frac{1}{2\pi}\int_{-\infty}^{\infty}
 \Im(\la)\ln\Big(\det\big(\mu^{\prime}(t)\big)\Big)\frac{dt}{|t-\la|^2}.
\end{align}
Formulas \eqref{As21}, \eqref{As34} and \eqref{As27} together with Proposition \ref{PnNond} yield the relations
\begin{align}\nn
\det\big(G_{\mu}(\la)^*G_{\mu}(\la)\big)&\geq \ov{\lim}_{r\to \infty}\det\big(2\pi \rho_{k_r}(\la,\ov{\la})\big)^{-1}
\\ &  \label{As35}
=\lim_{k\to \infty}\det\big(2\pi \rho_k(\la,\ov{\la})\big)^{-1}.
\end{align}
Using Lemma \ref{LaDet} and formulas \eqref{As8} (for $z=\la$) and \eqref{As35}, we derive \eqref{As9}.
\end{proof}
\begin{Rk}\label{RkEx}
The considerations of the Step 3 in the proof of Theorem \ref{TmMAs} $($based on the matrix balls results in Appendix \ref{MB}$)$
show that $\vp(z)$ satisfying \eqref{As7} always exists even without the requirements \eqref{S11} or \eqref{As8+}.
\end{Rk}
\begin{Rk}\label{RkRho}  According to Part II of Theorem \ref{TmMAs}, condition \eqref{As8+} $($for some $\la\in \BC_+)$
yields Szeg{\H o} inequality for certain $\vp(z)$ satisfying \eqref{As7}. According to Part I of Theorem \ref{TmMAs},
this implies \eqref{As8}. Therefore, if  only  \eqref{As8+} holds for some $\la\in \BC_+$, it also holds for all $\la\in \BC_+$.
\end{Rk}
\begin{Rk}\label{Rkcml} Similar to the special case of Toeplitz matrices in \cite{ALS-ConstrAppr}, asymptotics \eqref{As9} yields
asymptotics of 
$$\clm(k,z, \ov{\la}):=\mfa(S_k,z)J\mfa(S_k,\la)^*  \quad (z,\la\in \BC_+)$$ 
for the general $S$-nodes case.
\end{Rk}

Theorem \ref{TmMAs} is easily applied  to the case of Hankel matrices considered in Subsection \ref{Hank}
and related to the orthogonal polynomials on the real line (OPRL case).
\begin{Cy}\label{CyH} Let Hankel matrices $H(n)$ $(n\in \BN)$ of the form \eqref{c32}
be positive definite $($i.e., let the condition $H(n)>0$ hold$)$. Let the $S$-nodes $\{A_n=A(n), S_n=H(n),\Pi_n=\Pi(n)\}$
be given by \eqref{H2}--\eqref{H4} and assume that \eqref{As8+} is valid for some $\la \in \BC_+$.

 Then, the conditions of Theorem \ref{TmMAs}, Part I are fulfilled and relations
\eqref{As8} are valid for the matrix functions $\vp(z)$ satisfying \eqref{S11} and \eqref{As7}. 
Moreover, the equalities \eqref{As9} hold for each $\la \in \BC_+$ and the corresponding $\vp(z)$
satisfying \eqref{As7} $($and constructed in the proof of Theorem \ref{TmMAs}$)$.
\end{Cy}
\begin{proof}. It is easy to see that the matrices $I-zA_k$ are invertible for all complex $z$
and that the conditions of  Proposition \ref{PnPM} are fulfilled. 
The existence of $\vp(z)$
satisfying \eqref{As7} follows, for instance, from Remark \ref{RkEx}.  
For $\vp(z)$ satisfying \eqref{S11} as well, the conditions of Part I of Theorem \ref{TmMAs}, are fulfilled.

According to Remark \ref{RkRho}, \eqref{As8+} holds for all $\la\in \BC_+$. Hence, 
the conditions of Part II of Theorem \ref{TmMAs} are also fulfilled for all $\la\in \BC_+$.
\end{proof}
Note that the uniqueness of $\vp(z)$ satisfying \eqref{As7} immediately provides an essentially stronger result
than the one in the Theorem \ref{TmMAs}.
\begin{Cy}\label{CyUniq} Let the conditions of Proposition \ref{PnPM} be fulfilled for the $S$-nodes
$\{A_k,S_k,\Pi_k\}$ $(1 \leq k < \infty)$,  let
$\mfa(S_k,z)$ be meromorphic in $\BC_-$ and assume that the
inequality \eqref{As8+} holds for some
fixed $\la\in \BC_+$.
Let  the matrix function
$\vp(z)$ satisfying \eqref{As7} be unique. Then, for $\mu$ from the Herglotz representation of $\vp(z)$
and all $\la\in \BC_+$ 
we have the following asymptotics:
\begin{align}& \label{As37}
\lim\limits_{k\to \infty}\rho_k\big(\la, \ov{\la}\big)^{-1} = 2 \pi  G_{\mu}(\la)^* G_{\mu}(\la) .
\end{align} 
\end{Cy}
\begin{proof}. According to Remark \ref{RkRho}, \eqref{As8+} holds for all $\la \in \BC_+$.
Now, the uniqueness of $\vp(z)$ satisfying \eqref{As7}  and Part II of Theorem \ref{TmMAs}
imply that \eqref{As9} (or, equivalently\eqref{As37}) holds for $\mu$ from the Herglotz representation
of $\vp(z)$ and all $\la \in \BC_+$.
\end{proof}
\begin{Rk}\label{RkReal}
The asymptotics of $\rho_k(z,z)^{-1}$ $(z\in \BR)$ is of interest as well and is connected with the jump
of $\mu(t)$ at $t=z$ \cite{ALS87}.
\end{Rk}

\appendix

\section{On matrix ball representation}\label{MB}
\setcounter{equation}{0}
Representation of the linear-fractional transformations of the \eqref{S10} type in the form of matrix balls (or Weyl disks)
is a well-known and useful tool (see, e.g., \cite{Pot1} and some references therein).
For the purposes of self-sufficiency, we present it here under conditions considered in Section \ref{Ine}.
\begin{Pn}\label{PnMB}  Let a triple $\{A,S,\Pi\}$ form an $S$-node, let  relations \eqref{S9+}  hold, and let
the operators $I-zA$ have bounded inverses for $z$ in the domains  $\{z: \, \Im(z)\leq 0\}$ and $\{z: \, \Im(z)>0, \,\, |z| \geq r_0\}$ for some $r_0>0$.
Assume that $\mfa(S,z)=\{\mfa_{ij}(S,z)\}_{i,j=1}^2$ is given by \eqref{S9} and $\{R(z),Q(z)\}$ are  nonsingular pairs with property-$J$.

Then, the values of the matrix functions $\vp(z)$ obtained via linear-fractional transformations \eqref{S10}  form at each $z\in \BC_+$ $($such that $I-zA$ has a bounded inverse$)$
the following matrix ball:
\begin{align}& \label{B0}
\vp(z)=\I\big(- \rho(\ov{z},z)\ \big)^{-1}\al_{12}(z)-\big(- \rho(\ov{z},z)\ \big)^{-1/2}u\,\rho(z,\ov{z})^{-1/2},
\end{align}
where $u^*u\leq I_p$, $\I\big(- \rho(\ov{z},z)\ \big)^{-1}\al_{12}(z)$ is the so called centre of the ball, $\big(- \rho(\ov{z},z)\ \big)^{-1/2}$ is the
left radius and $\rho(z,\ov{z})^{-1/2}$ is the right radius.
\end{Pn}
\begin{proof}.  It follows from \eqref{S12} that 
\begin{align}& \label{B1}
\mfa(S,z)^{-1}=J \mfa(S,\ov{z})^*J.
\end{align}
We set
\begin{align}& \label{B2}
\al(S,z)=\{\al_{ik}(z)\}_{i,k=1}^2:=\big( \mfa(S,{z})^{-1}\big)^*J\mfa(S,{z})^{-1}.
\end{align}
Formula \eqref{S10} may be rewritten as
\begin{align}& \label{B3}
\mfa(S,{z})^{-1}\begin{bmatrix}-\I \vp(z) \\ I_p\end{bmatrix}=\begin{bmatrix}R(z) \\ Q(z)\end{bmatrix}\bigl( {\mathfrak A}_{21}(z ) R(z)+{\mathfrak A}_{22}(z)Q(z)\bigr)^{-1}.
\end{align}
Thus,  we may rewrite the second relation in \eqref{c16} as
\begin{align}& \label{B4}
\begin{bmatrix}\I \vp(z)^* & I_p\end{bmatrix}\al(z)\begin{bmatrix}-\I \vp(z) \\ I_p\end{bmatrix}\geq 0.
\end{align}
In view of \eqref{B1}, \eqref{B2} and \eqref{I8+}, we have
\begin{align}& \label{B5}
\al_{11}(z)=\big(\mfa(S,\ov{z})J\mfa(S,\ov{z})^*\big)_{22}=\rho(\ov{z},z).
\end{align}
Hence, \eqref{I8} implies that $\rho(\ov{z},z)<0$ for $z\in \BC_+$ and $\al_{11}(z)$ is invertible.
Therefore, $\al(z)$ may be represented in the form
\begin{align}\nn
\al(S,z)=&\begin{bmatrix}\al_{11}(z) \\ \al_{21}(z)\end{bmatrix}\al_{11}(z)^{-1}\begin{bmatrix}\al_{11}(z) & \al_{12}(z)\end{bmatrix}
\\ & \label{B6}
+\begin{bmatrix}0 & 0 \\ 0 & \al_{22}(z) - \al_{21}(z)\al_{11}(z)^{-1} \al_{12}(z)\end{bmatrix}.
\end{align}
It is easily calculated (see, e.g., \cite{SaL1+}) that
$$\al_{22}(z) - \al_{21}(z)\al_{11}(z)^{-1} \al_{12}(z)=\big(\big(\al(S,z)^{-1}\big)_{22}\big)^{-1},$$
where $\big(\al(S,z)^{-1}\big)_{22}$ is the  lower right $p\times p$ block of $\al(S,z)^{-1}$. 
We note that according to  \eqref{I8+} and \eqref{B2} we have
$\big(\al(S,z)^{-1}\big)_{22}=\rho(z,\ov{z})>0$. Hence, $\big(\al(S,z)^{-1}\big)_{22}$ in the formula
above is, indeed, invertible and we have
\begin{align}& \label{B7}
\al_{22}(z) - \al_{21}(z)\al_{11}(z)^{-1} \al_{12}(z)=\rho(z,\ov{z})^{-1}>0.
\end{align}
Using \eqref{B5}--\eqref{B7}, we rewrite \eqref{B4} in the form
\begin{align}&\label{B8}
\rho(z,\ov{z})^{-1}\geq   \big( \rho(\ov{z},z)\vp(z)+ \I\al_{12}(z)\big)^*\big(- \rho(\ov{z},z)\ \big)^{-1} \big( \rho(\ov{z},z)\vp(z)+ \I\al_{12}(z)\big).
\end{align}
Setting 
\begin{align}&\label{B9}
u(z):=\big(- \rho(\ov{z},z)\ \big)^{-1/2} \big( \rho(\ov{z},z)\vp(z)+ \I\al_{12}(z)\big)\rho(z,\ov{z})^{1/2},
\end{align}
we see that \eqref{B0} holds for this $u$. (In fact, \eqref{B0} is equivalent to \eqref{B9}.) Moreover, \eqref{B8} yields $u(z)^*u(z)\leq I_p$.
In other words, the matrices $\vp(z)$ given by \eqref{S10} belong to the matrix ball \eqref{B0}, where $u(z)^*u(z)\leq I_p$.

On the other hand, given any contractive matrix $u$, we see that $\vp(z)$ of the form
\eqref{B0} at $z$ is generated (via \eqref{S10}) by the nonsingular pair
\begin{align}& \label{B10}
\begin{bmatrix}R \\ Q\end{bmatrix}=\mfa(S,{z})^{-1}\begin{bmatrix}-\I \vp(z) \\ I_p\end{bmatrix}.
\end{align}
In order to show that this pair has property-$J$, we rewrite \eqref{B0} in the form \eqref{B9}
and derive \eqref{B8} from $u^*u\leq I_p$. Now, \eqref{B4} follows from \eqref{B8}.
In view of \eqref{B2} and \eqref{B4}, we obtain the  property-$J$ of the pair $\{R,Q\}$ given by
\eqref{B10}. That is, each matrix of our matrix ball coincides with the value of some $\vp\in \cln\big(\mfa(S)\big)$
at $z$.
\end{proof}
\section{On a certain limit inequality}\label{App}
\setcounter{equation}{0}
Passage to the limit under the integral sign is a classical topic in analysis (see, e.g., \cite[Ch. 6]{Nat}). 
In this appendix, we consider  matrix functions $\clp_k(t)$  such that $(1+t^2)^{-1}\clp_k(t)$ are integrable on $\BR$.
We will derive the following proposition, which is required in the proof of Theorem \ref{TmMAs}.
\begin{Pn} \label{LIn} Let a $p\times p$ nondecreasing matrix function $\mu(t)$ satisfy the inequality in \eqref{c18} and a sequence of 
$p\times p$ matrix functions $\clp_k(t)$ $(k\in \BN,$  $t \in \BR)$ satisfy relations
\begin{align} & \label{Ap1} 
\clp_k(t)\geq 0 , \quad \int_{-\infty}^{\infty}(1+t^2)^{-1}\clp_k(t)dt \leq M I_p \quad (k\in \BN, \,\, M >0);
\\ & \label{Ap2}
\lim_{k\to \infty}\int_{-\infty}^x (1+t^2)^{-1}\clp_k(t)dt=\gamma_{\clp}+\int_{-\infty}^x (1+t^2)^{-1}d\mu(t) \quad (\gamma_{\clp} \geq 0).
\end{align} 
Then, for all the continuous and bounded functions $f(t)> 0$ $(t \in \BR)$, for $a\in (\BR\cup\{-\infty\})$  
and for $b\in (\BR\cup\{\infty\})$ $(a<b)$,
we have
\begin{align}& \label{Ap3}
\displaystyle{\ov{\lim}_{k\to \infty}}\int_a^b f(t)\ln\Big(\det\big(\clp_k(t)\big)\Big)\frac{dt}{1+t^2} \leq \int_a^b f(t)\ln\Big(\det\big(\mu^{\prime}(t)\big)\Big)\frac{dt}{1+t^2}.
\end{align} 
\end{Pn} 
\begin{Rk} \label{RkLIn} The proposition above remains valid if $f(t)$ turns to $0$ at isolated points. We assume that \eqref{Ap3} holds if its left-hand side
is $-\infty$. It will follow from the proof of the proposition that the left-hand side of \eqref{Ap3} equals
$-\infty$ if the right-hand side of \eqref{Ap3} equals $-\infty$.
\end{Rk}
\begin{Rk}\label{Helly} The classical Helly's selection theorem $($scalar case$)$ is easily generalised to Helly's selection theorem $($matrix case$)$.
\end{Rk}
\begin{Rk}\label{HellyM} According to Helly's selection theorem  $($matrix case$)$, there is a subsequence of the sequence of integrals on  the left-hand side of
\eqref{Ap2}, which converges to some non-decreasing matrix function $\mu_0(t)$. Setting
\begin{align}& \label{mu1}
\mu(t):=\int_0^t(1+u^2)d\mu_0(u),
\end{align} 
we derive
\begin{align}& \label{mu2}
\int_{-\infty}^x (1+t^2)^{-1}d\mu(t) =\int_{-\infty}^x d\mu_0(t)=\mu_0(t)-\g_{\clp}, 
\end{align} 
where $\g_{\clp}:=\mu_0(-\infty) \geq 0$. Therefore, the convergence in \eqref{Ap2} $($under conditions \eqref{Ap1}$)$ is fulfilled for some subsequence at least.
\end{Rk}
In order to prove Proposition \ref{LIn}, we need two lemmas.
\begin{La}\label{Lapn1}
Let the conditions of Proposition \ref{LIn} hold and assume that $p=1$. Then, for all the continuous and bounded functions $f(t)\geq 0$ $(t \in \BR)$, for $a\in (\BR\cup\{-\infty\})$  
and for $b\in (\BR\cup\{\infty\})$ $(a<b)$,
we have
\begin{align}& \label{Ap4}
\displaystyle{\ov{\lim}_{k\to \infty}}\int_a^b f(t)\ln\big(\clp_k(t)\big)\frac{dt}{1+t^2} \leq \int_a^b f(t)\ln\big(\mu^{\prime}(t)\big)\frac{dt}{1+t^2}.
\end{align} 
\end{La}
\begin{proof}. We set
\begin{equation} \label{Ap5}
\clp_{kn}(t)=\left\{ \begin{array}{l} \clp_k(t) \quad {\mathrm{for}} \quad \ve_n <\clp_k(t)<M_n, \\ \ve_n \quad {\mathrm{for}} \quad \clp_k(t)\leq \ve_n,
\\ M_n \quad {\mathrm{for}} \quad \clp_k(t) \geq M_n, \end{array}\right.  
\end{equation} 
where $n\in \BN$, $0<\ve_n<1, \,\, M_n>1,$ and $\ve_n\to 0, \,\, M_n\to \infty$ for $n\to \infty$. We also  introduce the function $\xi$:
\begin{align}& \label{Ap6}
\xi(t)=\int_{-\infty}^t (1+u^2)^{-1}du, \quad d\xi =(1+t^2)^{-1}dt.
\end{align}
Clearly, the inverse function $t(\xi)$ (or $t(v)$) exists for $\xi\in [0,\pi)$ and the functions in the sequences
\begin{align}& \label{Ap7}
\int_0^\xi \clp_{kn}\big(t(v)\big)dv \quad (k\in \BN)\quad {\mathrm{and}} \quad \int_0^\xi \ln\Big(\clp_{kn}\big(t(v)\big)\Big)dv  \quad (k\in \BN)
\end{align} 
are bounded and absolutely equicontinuous in the terminology of \cite[ p.~20]{Priv}. According to an analogue 
of Arzel\`a--Ascoli theorem from \cite[p. 21]{Priv}, there are subsequences of these sequences, which uniformly converge
on $[0,\pi)$ to absolutely continuous functions. Without loss of generality, we may consider instead sequences \eqref{Ap7},
which uniformly converge to integrals of the functions ${}^1\clp_n$ and $\ln({}^2\clp_n)$, respectively, that is, to
\begin{align}& \label{Ap8}
\int_0^\xi {}^1\clp_n\big(t(v)\big)dv=\int_{-\infty}^t{}^1\clp_n(u)\frac{du}{1+u^2} \quad {\mathrm{and}} \\ & \label{Ap9}
 \int_0^\xi \ln\Big({}^2\clp_n\big(t(v)\big)\Big)dv =\int_{-\infty}^t\ln\big({}^2\clp_n(u)\big)\frac{du}{1+u^2} \, ,
\end{align} 
and such that the expressions for the corresponding $\clp_k$ on the left-hand side of  \eqref{Ap4} tend to the upper limit .
 Taking into account, which sequences converge to the integrals \eqref{Ap8} and \eqref{Ap9}, and  the formula on arithmetic and geometric means in integral form
 (see, e.g., \cite[(6.7.5)]{HLP}), we have
\begin{align}& \label{Ap10}
\frac{1}{b_1-a_1}\int_{a_1}^{b_1}\ln\Big({}^2\clp_n\big(t(v)\big)\Big)dv \leq \ln\left(\frac{1}{b_1-a_1}\int_{a_1}^{b_1}{}^1\clp_n\big(t(v)\big)dv\right)
\end{align}
for $0\leq a_1<b_1<\pi$. The inequalities \eqref{Ap10} yield 
\begin{align}& \label{Ap11}
\ln \big({}^2\clp_n(t)\big)\leq \ln \big({}^1\clp_n(t)\big) \quad (-\infty<t<\infty).
\end{align}
It follows from \eqref{Ap11} that
\begin{align}& \label{Ap12}
\displaystyle{\ov{\lim}_{k\to \infty}}\int_a^b f(t)\ln\big(\clp_{kn}(t)\big)\frac{dt}{1+t^2} \leq \int_a^b f(t)\ln\Big({}^1\clp_n(t)\Big)\frac{dt}{1+t^2}.
\end{align} 

Since $\big(x^{- \vk}\ln(x)\big)^{\prime}=x^{-1- \vk}\big(1-\vk \ln(x)\big)$, we have  $\ln(x)\leq x^{\vk}\big/ (\vk \E)$ for $\vk>0$ and $x\in [1,\infty)$.
Setting $\vk=1$ and using \eqref{Ap1}, we derive
\begin{align}& \label{Ap13}
\int_a^{b}\ln_{+}\big(\clp_k(t)\big/\clp_{kn}(t)\big)\frac{dt}{1+t^2} \leq \frac{1}{\E M_n}\int_{a}^{b}\frac{\clp_k(t)dt}{1+t^2}\leq M\big/(\E M_n),
\end{align}
where $\ln_+ (a)= 0$ for $0<a\leq 1$ and $\ln_+ (a)= \ln(a)$ for $a > 1$. Put
\begin{equation} \label{Ap5'}
\wt \clp_{kn}(t)=\left\{ \begin{array}{l} \clp_k(t) \quad {\mathrm{for}} \quad  \clp_k(t)>\ve_n, \\ \ve_n \quad {\mathrm{for}} \quad \clp_k(t)\leq \ve_n.
 \end{array}\right. 
\end{equation} 
According to Helly's selection theorem, there are subsequences (and without a loss of generality we again consider them as a
sequence) such that $\break \int_{-\infty}^t(1+u^2)^{-1}\wt \clp_{kn}(u)du$ converges to $\g+\int_{-\infty}^t(1+u^2)^{-1}d\mu_n(u)$,
where $\mu_n(t)$ is nondecreasing. We set $\wt \clp_{n}(t)=\mu_n^{\prime}(t)$, where $\mu_n^{\prime}$ is the derivative of the absolutely
continuous part of $\mu_n$.

Using \eqref{Ap13} and \eqref{Ap5'}, we obtain
\begin{align}& \nn
\uv{\lim}_{n\to \infty}{\ov{\lim}_{k\to \infty}}\int_a^b f(t)\ln\big(\clp_{kn}(t)\big)\frac{dt}{1+t^2} 
\\ & \label{Ap14}
=\uv{\lim}_{n\to \infty}{\ov{\lim}_{k\to \infty}}\int_a^b f(t)\ln\big(\wt \clp_{kn}(t)\big)\frac{dt}{1+t^2} 
\\ & \nn
\geq {\ov{\lim}_{k\to \infty}}\int_a^b f(t)\ln\big( \clp_{k}(t)\big)\frac{dt}{1+t^2} .
\end{align} 
It is easy to see that $\wt \clp_{n}(t) \geq {}^1\clp_n(t)$, and so \eqref{Ap12} yields
\begin{align}& \nn
\uv{\lim}_{n\to \infty}{\ov{\lim}_{k\to \infty}}\int_a^b f(t)\ln\big(\clp_{kn}(t)\big)\frac{dt}{1+t^2} 
\leq {\uv{\lim}_{n\to \infty}}\int_a^b f(t)\ln\big(\wt \clp_{n}(t)\big)\frac{dt}{1+t^2} .
\end{align} 
Hence, taking into account \eqref{Ap14} we derive
\begin{align}& \label{Ap16}
{\ov{\lim}_{k\to \infty}}\int_a^b f(t)\ln\big( \clp_{k}(t)\big)\frac{dt}{1+t^2} 
\leq {\uv{\lim}_{n\to \infty}}\int_a^b f(t)\ln\big(\wt \clp_{n}(t)\big)\frac{dt}{1+t^2} .
\end{align} 
Since $\wt \clp_{n}(t) \leq \mu^{\prime}(t)+\ve_n$, formula \eqref{Ap16} implies \eqref{Ap4}.
\end{proof}

In order to study the case $p>1$ we will need the following inequality for 
the  $p\times p $ matrices $B_k\geq 0$ $(1\leq k\leq m)$:
\begin{align}& \label{Ap17}
\sum_{k=1}^m\big(\det B_k\big)^{1/p} \leq \left(\det\Big(\sum_{k=1}^m B_k\Big)\right)^{1/p}.
\end{align}
Indeed, \eqref{Ap17} easily follows (by induction) from Minkowski  inequality 
\begin{align}& \label{Ap21}
\big(\det B_1\big)^{1/p}+\big(\det B_2\big)^{1/p}\leq \big(\det(B_1+ B_2)\big)^{1/p}.
\end{align}
(For the Minkowski  inequality  see, e.g.,  \cite{HLP}.)

Taking into account \eqref{Ap17}, we will prove the lemma below.
\begin{La}\label{La5'} Let a sequence of uniformly bounded $p\times p$ matrix functions $\clp_k(t)\geq 0$
satisfy \eqref{Ap2} and assume that
\begin{align}  & \label{Apx1}
\lim_{k\to \infty}\int_{-\infty}^x (1+t^2)^{-1}\big(\det \clp_k(t)\big)^{1/p}dt=\int_{-\infty}^x (1+t^2)^{-1}\tau(t)dt 
\end{align} 
for $-\infty<x<\infty, \,\,\tau(t) \geq 0$. Then, we have
\begin{align}& \label{Apx2}
\tau(t)\leq \big(\det \mu^{\prime}(t)\big)^{1/p}.
\end{align}
\end{La}
\begin{Rk}\label{RkSeq} The convergences \eqref{Ap2} and \eqref{Apx1}, with certain absolutely continuous $\mu$ $(\mu^{\prime}(t)\geq 0)$ and $\tau(t) \geq 0$,
follow for some subsequence of  a uniformly bounded sequence of $p\times p$ matrix functions $\clp_k(t)\geq 0$  from the considerations 
at the beginning of the proof of Lemma \ref{Lapn1} in any case.
\end{Rk}
\begin{proof} of Lemma \ref{La5'}. Assume that \eqref{Apx2} does not hold. Then, there is some set $\clu$ of  nonzero measure $\clm(\clu)=\int_{\clu}\frac{dt}{1+t^2}$
such that the sequence $\clp_k$ tends there to $\mu^{\prime}(t)$, $\big(\det \clp_k(t)\big)^{1/p}$ tends to $\tau(t)$ (as in the conditions of the lemma) and
\begin{align}& \label{Apx3}
\tau(t)\geq \big(\det \mu^{\prime}(t)\big)^{1/p}+\ve_1, \quad  \ve_2 I_p\leq \mu^{\prime}(t)\leq \clm_2 I_p \quad (\ve_1,\ve_2>0).
\end{align}
We  will require (for the choice of $\clu$) additionally  that
\begin{align}& \label{Apx4}
-\de I_p \leq \mu^{\prime}(t_0)-\mu^{\prime}(t) \leq \de I_p
\end{align}
for some $t_0\in \clu$, almost all $t\in \clu$ and sufficiently small values $\de>0$. (It may be done, although this additional requirement
reduces the initial set.)
We choose $\de <\ve_2/2$ so that
for $\clp:\equiv \mu^{\prime}(t_0)-\de I_p$ we have
\begin{align}& \label{Apx5}
\mu^{\prime}(t)>2\de I_p, \,\, \mu^{\prime}(t) \geq \clp>0 \,\, (t\in \clu), \quad \int_{\clu}\big(\mu^{\prime}(t) - \clp\big)\frac{dt}{1+t^2}\leq 2\de \clm(\clu)I_p.
\end{align}
We set 
$$B:= \int_{\clu}\mu^{\prime}(t) \frac{dt}{1+t^2},$$ 
and obtain from \eqref{Apx5} that
\begin{align}& \label{Apx6}
\big(\det(B- 2\de \clm(\clu)I_p)\big)^{1/p}\leq \int_{\clu}(\det \clp)^{1/p}\frac{dt}{1+t^2}\leq \int_{\clu}\big(\det \mu^{\prime}(t)\big)^{1/p}\frac{dt}{1+t^2}.
\end{align}
Approximating $\clp_k(t)$ by matrices taking finite numbers of values and using \eqref{Ap17}, we derive
\begin{align}& \label{Apx7}
\int_{\clu}\big(\det \clp_k(t)\big)^{1/p}\frac{dt}{1+t^2}\leq \left(\det \int_{\clu} \clp_k(t)\frac{dt}{1+t^2}\right)^{1/p}.
\end{align}
For each $\de_0>0$, there is $k_0$ such that
\begin{align}& \label{Apx8}
\int_{\clu}\clp_k(t)\frac{dt}{1+t^2}\leq B+\de_0 I_p \quad {\mathrm{for}} \quad k>k_0.
\end{align}
Recall that $\big(\det \clp_k(t)\big)^{1/p}$ converges to $\tau(t)$. Hence, relations \eqref{Apx7} and \eqref{Apx8} imply that
\begin{align}& \label{Apx9}
\int_{\clu}\tau(t)\frac{dt}{1+t^2}\leq (\det B)^{1/p}.
\end{align}
It follows from \eqref{Apx6} and \eqref{Apx9} that
\begin{align}& \label{Apx10}
\int_{\clu}\tau(t)\frac{dt}{1+t^2}\leq \int_{\clu}\big(\det \mu^{\prime}(t)\big)^{1/p}\frac{dt}{1+t^2}+O(\de)\clm(\clu),
\end{align}
where, as usually, $O(\de)\leq \clm_3 \de$ for some $\clm_3>0$ and sufficiently small $\de$. Thus, \eqref{Apx10}
contradicts the first inequality in \eqref{Apx3}, and \eqref{Apx2} is proved by negation.
\end{proof}
Now, we will prove Proposition \ref{LIn} using Lemmas \ref{Lapn1} and \ref{La5'}.
\begin{proof} of Proposition \ref{LIn}.
Let us represent $\clp_k$ in the form 
\begin{align}& \label{Ap22}
\clp_k(t)=U_k(t)\cld_k(t)U_k(t), \quad \cld_k(t)=\diag\{d_1(t,k),\ldots , d_p(t,k)\},
\end{align}
where $U_k(t)$ are unitary, $\diag$ stands for a diagonal matrix and $d_{\ell}(t,k)$ are the entries
of $\cld_k(t)$ on the main diagonal. We set
\begin{align}& \label{Ap23}
\clp_{kn}(t)=U_k(t)\cld_{kn}(t)U_k(t);
\\ & \label{Ap23+} 
 \cld_{kn}(t)=\diag\{d_1(t,k,n),\ldots , d_p(t,k,n)\},
\\ & \label{Ap24} 
d_{\ell}(t,k, n)= d_{\ell}(t,k) \quad {\mathrm{for}} \quad d_{\ell}(t,k)<M_n, 
\\ & \label{Ap25} 
 d_{\ell}(t,k, n)=M_n \quad {\mathrm{for}} \quad d_{\ell}(t,k) \geq M_n,
\end{align}
where the sequence $M_n$ increases to $\infty$. According to Remark \ref{RkSeq},
we may choose a subsequence $\clp_{k_r}$ such that the upper limit in \eqref{Ap3} 
is achieved for this subsequence and the conditions of Lemma \ref{La5'}
hold for  the sequences $\clp_{k_r n}$ (and some $\mu_{n}$, $\tau_n$). Without loss of generality, we may exclude the case, where
 $\det (D_{k_r}(x))=0$ on a set with nonzero measure because in the case of an infinite number of  $\clp_{k_r}$ with this
 property
the inequality \eqref{Ap3} is immediate (see Remark~\ref{RkLIn}).
Moreover, without loss of generality, we assume that $\clp_{k_r}$ is our sequence
$\clp_k$. Then,  the sequences $\big(\det \clp_{kn}(t)\big)^{1/p}$ satisfy conditions of Lemma \ref{Lapn1}, where $\tau_n$
stands in place of $\mu^{\prime}$. It follows that
\begin{align}& \nn
\displaystyle{\ov{\lim}_{k\to \infty}}\int_a^b f(t)\ln\big(\big(\det \clp_{kn}(t)\big)^{1/p}\big)\frac{dt}{1+t^2} \leq \int_a^b f(t)\ln\big(\tau_n(t)\big)\frac{dt}{1+t^2}.
\end{align}
Hence, taking into account Lemma \ref{La5'}, we derive
\begin{align}& \nn
\displaystyle{\ov{\lim}_{k\to \infty}}\int_a^b f(t)\ln\big(\det \clp_{kn}(t)\big)\frac{dt}{1+t^2} \leq \int_a^b f(t)\ln\big(\det\mu_n^{\prime}(t)\big)\frac{dt}{1+t^2}.
\end{align}
It is easy to see that $\mu_n^{\prime}(t) \leq \mu^{\prime}(t)$. Thus, the formula above yields
\begin{align}& \label{Ap26}
\displaystyle{\ov{\lim}_{k\to \infty}}\int_a^b f(t)\ln\big(\det \clp_{kn}(t)\big)\frac{dt}{1+t^2} \leq \int_a^b f(t)\ln\big(\det\mu^{\prime}(t)\big)\frac{dt}{1+t^2}.
\end{align}
Since $\det \big(D_k(t)\big)\not=0$, we have $\det \big(D_{kn}(t)\big)\not=0$ almost everywhere as well, that is, $D_{kn}$ is invertible. We set
\begin{align}& \label{Ap27}
\clp_{kn}^r(t):=U_k(t)D_k(t)D_{kn}(t)^{-1}U_k(t)^*.
\end{align}
Clearly, $\ln\big(\det \clp_{k}(t)\big)=\ln\big(\det \clp_{kn}(t)\big)+\ln\big(\det \clp_{kn}^r(t)\big)$. Hence, we have
\begin{align}\nn
\displaystyle{\ov{\lim}_{k\to \infty}}\int_a^b f(t)\ln\big(\det \clp_{k}(t)\big)\frac{dt}{1+t^2}\leq &\displaystyle{\ov{\lim}_{k\to \infty}}\int_a^b f(t)\ln\big(\det \clp_{kn}(t)\big)\frac{dt}{1+t^2}
\\ & \label{Ap28}
+M \int_a^b\ln\big(\det \clp_{kn}^r(t)\big)\frac{dt}{1+t^2},
\end{align}
where $M$ is the maximum of $f(t)$ on the interval of integration. From \eqref{Ap22}--\eqref{Ap25} and \eqref{Ap27}, we obtain
\begin{align}& \label{Ap29}
I_p \leq \clp_{kn}^r(t) \leq I_p+ \clp_{k}(t) / M_n.
\end{align}
Moreover, using Lemma 6.3 in \cite[Ch. II]{Roz} (after the change of variables as in \eqref{Ap6}, see also \eqref{Ap7}--\eqref{Ap9}), we derive the inequality
\begin{align}\nn
\int_a^b \ln\big(\det \clp_{kn}^r(t)\big)\frac{dt}{1+t^2}&=\int_{\xi_1}^{\xi_2} \ln\big(\det \clp_{kn}^r(t(v))\big){dv}
\\ & \label{Ap30}
\leq (\xi_2-\xi_1)\ln\left(\det \frac{1}{\xi_2-\xi_1}\int_{\xi_1}^{\xi_2}\clp_{kn}^r(t(v))\big)dv\right).
\end{align}
Taking into account \eqref{Ap28}--\eqref{Ap30}, we derive
\begin{align}& \nn
\displaystyle{\ov{\lim}_{k\to \infty}}\int_a^b f(t)\ln\big(\det \clp_{k}(t)\big)\frac{dt}{1+t^2}
\\ & \label{Ap31}
\leq \ov{\lim}_{n\to \infty}\displaystyle{\ov{\lim}_{k\to \infty}}\int_a^b f(t)\ln\big(\det \clp_{kn}(t)\big)\frac{dt}{1+t^2}.
\end{align}
Finally, the inequalities \eqref{Ap26} and \eqref{Ap31} yield \eqref{Ap3}.
\end{proof}
\section{Determinant inequality}\label{Det}
\setcounter{equation}{0}
Here,  we obtain the following lemma, which seems almost evident  but needs some proof in spite of that.
\begin{La}\label{LaDet} Let $\cla$ and $\clb$ be two $p\times p$ matrices such that
\begin{align}& \label{Ac1}
\cla>0, \quad \clb\geq 0, \quad \clb\not=0.
\end{align}
Then, the strict inequality $\det(\cla+\clb)>\det \cla$ is valid.
\end{La}
\begin{proof}. We denote the eigenvalues of $\cla$ by $z_k$ and the eigenvalues of $\cla+\clb$ by $\la_k$
and assume that $z_1\geq z_2\geq \ldots \geq z_p$ and $\la_1\geq \la_2\geq\ldots \geq \la_p$.
The eigenvalue $z_k$ may be written down in the form
\begin{align}& \label{Ac2}
z_k=\max_{V_k}\,\, \min_{h\in V_k,  \,\, h^*h=1}  \,\, h^*\cla h,
\end{align}
where $V_k$ are the $k$-dimensional subsets of  $\BC^p$ (see, e.g., \cite[p. 545]{Gant}).
The inequalities $\la_k\geq z_k$ and $\det(\cla+\clb)\geq\det \cla$ are immediate from \eqref{Ac1} and \eqref{Ac2}.

In order to derive the strict inequality for the determinants above, we use a unitary similarity transformation $U$,
which transforms $\cla$ into the diagonal matrix $\cla_{tr}:=\diag\{z_1,\ldots, z_p\}$ and transforms $\clb$ into the block matrix
\begin{align}& \label{Ac3}
\clb_{tr}=\begin{bmatrix} \cld & \clb_{12} \\ \clb_{21} & \br \clb\end{bmatrix}, \quad \cld=\diag\{d_1, \ldots, d_s\} \quad (d_1 \geq \ldots \geq d_s),
\end{align}
where $s$ is the number of the maximal eigenvalues of $\cla$ so that $\cla_{tr}$ consists of two diagonal blocks:
\begin{align}& \label{Ac4}
\cla_{tr}=\diag\{z_1I_s, \br \cla\}, \quad \br \cla:=\diag\{z_{s+1}, \ldots, z_p\}.
\end{align}
In view of \eqref{Ac2}--\eqref{Ac4}, the inequality $\cld\not=0$ yields $\la_1>z_1$ and so $\break \det(\cla+\clb)=\det(\cla_{tr}+\clb_{tr})>\det \cla$.
If $\cld=0$, we take into account  the inequality $\clb_{tr}\geq 0$ and derive $B_{12}=B_{21}^*=0$.
Therefore, the inequality $\det(\cla+\clb)>\det \cla$ is equivalent to the inequality $\det(\br \cla+\br \clb)>\det \br \cla$,
and we come to the matrices $\br \cla$ and $\br \clb$ (of the reduced order $p-s$), which may be treated in the same way
as $\cla$ and $\clb$.
Since $\clb \not=0$, at some step we will arrive at the case $\cld\not= 0$.
\end{proof}


\end{document}